\documentclass[10pt,twocolumn,twoside,nonpeerreview]{IEEEtran}
\usepackage{ifpdf}

\usepackage{cite}

\ifCLASSINFOpdf
  \usepackage[pdftex]{graphicx}
\else
\fi
\usepackage{amsmath}
\usepackage{array}
\usepackage{url}

\usepackage[square,sort,comma,numbers]{natbib}
\usepackage{bbm}
\usepackage{amssymb}
\usepackage{amsthm}
\usepackage{dsfont}
\usepackage[hidelinks]{hyperref}
\hypersetup{%
pdftitle={A Stochastic Resource-Sharing Network for Electric Vehicle Charging},
pdfauthor={Aveklouris, Vlasiou, Zwart},
pdfkeywords={Electric vehicle charging, distribution network, AC power flow model, linearized (simplified) Distflow, queueing theory, stochastic processes, fluid approximation},
bookmarks=true,%
bookmarksnumbered=true,
pdfstartview={FitH},
breaklinks=true,%
}%
\newcommand{\ind}[1]{\mathds{1}_{\{#1\}}}   
\newcommand{\E}[1]{\mathbb{E}[#1]}        

\newcommand{\Prob}{\mathbb{P}}

\newcommand{\blt}{\boldsymbol}
\newcommand{\R}{\mathds{R}}
\newcommand{\SR}{\overline{R}}

\newcommand{\ep}{\epsilon}

\DeclareMathOperator{\Node}{\mathcal{I}}
\DeclareMathOperator{\Edge}{\mathcal{E}}
\DeclareMathOperator{\Type}{\mathcal{J}}

\graphicspath{{../Figures/}}
\newtheorem{defn}{Definition}
\newtheorem{rem}{Remark}
\newtheorem{thm}{Theorem}
\newtheorem{prop}{Proposition}

\usepackage{color}
\newcommand{\angb}[1]{\textcolor{black}{#1}}

\begin{document}
\title{A Stochastic Resource-Sharing Network for Electric Vehicle Charging}
\author{Angelos~Aveklouris,~
        Maria~Vlasiou,~
        and~Bert~Zwart,~\IEEEmembership{Member,~IEEE}
\thanks{The research of
Angelos Aveklouris is funded by a TOP grant of the Netherlands Organization for
Scientific Research (NWO) through project 613.001.301. The research of Maria
Vlasiou is supported by the NWO MEERVOUD grant 632.003.002.
The research of Bert Zwart is partly supported by the NWO VICI grant 639.033.413.}
\thanks{A.\ Aveklouris (a.aveklouris@tue.nl) is with the Department
of Mathematics and Computer Science, Eindhoven University of Technology, The Netherlands.}
\thanks{M.\ Vlasiou (m.vlasiou@tue.nl) and B. Zwart (bert.zwart@cwi.nl) are with the Centrum Wiskunde en Informatica and the Department of Mathematics and Computer Science, Eindhoven University of Technology, The Netherlands.}}%
%
\markboth{ }%
{A Stochastic Resource-Sharing Network for Electric Vehicle Charging}
%



\maketitle
\begin{abstract}
We consider a  distribution grid used to charge electric vehicles such that voltage drops stay bounded. We model this as a class of resource-sharing networks, known as \textit{bandwidth-sharing networks} in the communication network literature.
We focus on resource-sharing networks that are driven by a class of greedy control rules that can be implemented in a decentralized fashion. For a large number of such control rules, we can characterize the performance of the system by a fluid approximation. This leads to a set of dynamic equations that take into account the stochastic behavior of EVs. We show that the invariant point of these equations is unique and can be computed by solving a specific ACOPF problem, which admits an exact convex relaxation. \angb{We illustrate our findings with a case study using the SCE 47-bus network and several special cases that allow for explicit computations.}
\end{abstract}
%
\IEEEpeerreviewmaketitle
\section{Introduction}
\IEEEPARstart{T}{he} rise of electric vehicles (EVs) is unstoppable due to factors such as the decreasing cost of batteries and various policy decisions \cite{Energy2017}. These vehicles need to be charged and will therefore cause congestion in distribution grids in the (very near) future.
\angb{Current bottlenecks are the ability to charge a battery at a fast rate and the number of charging stations, but this bottleneck is expected to move towards the current grid infrastructure. This is illustrated in \cite{2025scenario}, where the authors evaluate the impact of the energy transition on a real distribution grid in a field study, based on a scenario for the year 2025. The authors confront a local low-voltage grid with electrical vehicles and ovens and show that charging a small number of EVs is enough to burn a fuse. Additional evidence of congestion is reported in \cite{clement2010impact}.
}

This paper proposes to model and analyze such congestion by the use of a class of resource-sharing networks, which in the queueing network community are known as \textit{bandwidth-sharing networks}. Bandwidth-sharing networks are a specific class of queueing networks, where customers (in our context: EVs) need to be served simultaneously. Their service requires the concurrent usage of multiple ``servers'' (in our case: all upstream lines between the location of the EV and the feeder of the distribution grid). Determining how fast to charge each EV taking into account network stability and the randomness of future arriving EVs is one of the key problems in the analysis of distribution grids, leading to challenging mathematical problems.

Similar questions have been successfully answered in communication networks, where the set of feasible schedules is determined by the maximum amount of data a communication channel can carry per time unit, leading to the powerful concept of bandwidth-sharing networks \cite{massoulie1999bandwidth}. Bandwidth-sharing networks couple the important fields of network utility maximization with stochastic process dynamics \cite{kelly1997charging, kelly2014stochastic}. Apart from yielding various qualitative insights, they have been instrumental in the comparison and optimization of various data network protocols and even in the design of new protocols \cite{yi2008stochastic}. The stochastic analysis of bandwidth-sharing networks was initially restricted to some specific examples \cite{BonaldProutiere2002, Bonald2006}. By now,
fluid and diffusion approximations are available, which are computationally tractable \cite{KangKellyLeeWilliams2008, YeYao2012, BorstEgorovaZwart2014, remerova14, ReedZwart2014, VlasiouZhangZwart2016} and hold for a large class of networks.

From a mathematical viewpoint, the present paper is influenced by \cite{remerova14}. However, in the setting of charging EVs, an important constraint that needs to be satisfied is to keep the voltage drops bounded, making the bandwidth-sharing network proposed in this paper different from the above-mentioned works. This also causes new technical issues, as the capacity set can be non-polyhedral or even non-convex. The first paper to suggest the class of bandwidth-sharing networks in the content of EV charging is \cite{carvalho2015critical}, where simulation studies were conducted to explore stability properties.

\angb{
Though simulation studies have value, they can be time-consuming and do not always offer qualitative insights. Other computational techniques, such as the numerical solution of Markov chains, quickly require the computation of trillions of convex programs - even for networks of less than 10 nodes and are therefore not a practical tool to assess the performance of distribution networks of realistic systems.}

\angb{
Our work is focusing on the development of a more tractable fluid model, using state-of-the-art techniques from stochastic network theory, in a model setting which is a significant extension of \cite{carvalho2015critical}: we allow for load limits, finitely many parking spaces, deadlines (associated with parking times), and do not make any assumption on the joint distribution of the parking time and the demand for electricity. More importantly, despite our assumptions leading to an intricate class of infinite-dimensional Markov processes (also known as measure-valued processes), we obtain a number of mathematical results that are computationally tractable and in some cases explicit. Our work can also be seen as an extension (to the case of stochastic dynamics) of the recent works \cite{ardakanian2013, fan2012}. In those papers it is shown, using arguments similar to the seminal work \cite{kelly1997charging}, how algorithms like proportional fairness emerge in decentralized EV charging. Our class of controls contain proportional fairness as a special case. Thus, our fluid model approach can serve as an additional exploratory tool for planning purposes.}

We develop a fluid approximation for the number of uncharged EVs (for the single-node Markovian case see \cite{aveklouriselectric, aveklouris2018bounds}), allowing the dynamics of the stochastic model to be approximated with a deterministic model. This model is still quite rich, as it depends on the joint distribution of the charging requirements and parking times. We show that the invariant point of this dynamical system is unique and can be characterized in a computationally friendly manner by formulating a nontrivial AC optimal-power-flow problem (ACOPF), which is tractable as its convex relaxation is exact in many cases. When we replace the AC load flow model with the simpler linearized Distflow \cite{Baran89}, we obtain more explicit results, as the capacity set becomes polyhedral. \angb{We illustrate our main result in a case study using the SCE 47-bus network in order to study the performance of the network. We present results for both the aggregated system and the individual nodes.}

\angb{In Appendix~\ref{Additional numerical examples},
we compare the fluid approximation with its original stochastic model and also compare across load flow models.
Our computational study, though certainly not exhaustive, indicates that the relative errors between the fluid model and the stochastic model are less than 10\% in most cases.}

For the class of weighted proportional fairness controls where the weights are chosen as function of the line resistances, we derive a fairness property. In this case, all EVs are charged at the same rate, independent of the location in the network, while keeping voltage drops bounded. When the weights are instead chosen equally, we can even derive an explicit formula for the invariant distribution of the original stochastic process. Specifically, we show that under certain assumptions, our network behaves like a multiclass processor-sharing queue. Such properties have proven quite fruitful in other areas of engineering, particularly in the analysis of computer systems \cite{Kleinrock1976}, communication networks \cite{BonaldProutiere2002}, and wireless networks \cite{borst05}.

EVs can be charged in several ways. Our setup can be seen as an example of slow charging, in which drivers typically park their EV and are not physically present during charging (but are busy shopping, working, sleeping, etc). For \emph{queueing} models focusing on fast charging, we refer to \cite{bayram2013electric, yudovina2015socially}. Both papers consider
a gradient scheduler to control delays, but do not consider physical load flow models as is\break done here. \cite{Swapping17} presents a queueing model for battery swapping. \cite{Turitsyn2010} is an early paper on a queueing analysis of EV charging, focusing on designing safe (in term of voltage drops) control rules with minimal communication overhead.

We can only provide a small additional (i.e., non-queueing) sample of the already vast but still emerging literature on EV charging. The focus of \cite{su2012performance} is on a specific parking lot and presents an algorithm for optimally managing a large number of plug-in EVs. Algorithms to minimize the impact of plug-in EV charging on the distribution grid are proposed in \cite{sortomme2011coordinated}.
In \cite{li2012modeling}, the overall charging demand of plug-in EVs is considered. Mathematical models where vehicles communicate beforehand with the grid to convey information about their charging status are studied in \cite{said2013queuing}. In \cite{kempker2016optimization}, EVs are the central object and a dynamic program is formulated that prescribes how EVs should charge their battery using price signals.

The present paper aims to illustrate how state-of-the-art methods from the applied probability and queueing communities can contribute to the analysis of the interplay between EV charging and the analysis of congestion in distribution networks. Our analysis does not take into account other important features in distribution networks touched upon in some of the above-cited works. In particular, it would be useful to incorporate smart appliances/buildings/meters, rooftop solar panels, and other sources of electricity demand. We think that the tractability and generality of our formulation and the established connection with an OPF problem makes our framework promising towards a comprehensive stochastic network model of a distribution grid.

\angb{
The paper is organized as follows. In Section \ref{Sec:model description}, we provide a detailed model description --- in particular we introduce our stochastic model, the class of charging controls, and the load flow models.
In Section~\ref{sec:FA}, we consider our model in full generality. We present a fluid model of our system and show that the associated dynamic equations have a unique invariant point, which is shown to be stable under an additional assumption. A case study using the SCE 47-bus network is presented in Section~\ref{sec:case Study}.
In Section~\ref{Additional examples under a line topology}, we show explicit  results under additional assumptions and we apply our results to design a control rule that maximizes the fraction of successful charges.}

\section{Model description}\label{Sec:model description}
In this section, we provide a detailed formulation of our model and explain various notational conventions that are used in the remainder of this work.
A notational convention is that all vectors and matrices are denoted by bold letters. \angb{We present the basic nomenclature in Table~\ref{table}.}

\begin{table}[!h]
\caption{Basic notation}
\label{table}
\centering
\resizebox{\columnwidth}{!}{%
\begin{tabular}{|c c|| c c|}
  \hline
$I$ & \# nodes  &$J$& \# EV types\\  \hline
$\epsilon_{ik}$ & edge $(i,k)$  &$\Node(i)$& subtree rooted at node $i$\\  \hline
$\mathcal{P}(i)$ & path from node $i$ to feeder  &$K$& \# parking spaces at nodes\\  \hline
$\lambda$ & arrival rate of Poisson process &$\gamma$ &  modified arrival rate due to finite parking spaces\\  \hline
$\blt{Q}$ & \# of total EVs at nodes &$\blt{Z}$& \# of uncharged EVs\\  \hline
$u(\cdot)$ & utility function &$\blt{z}^*$& fluid approximation of uncharged EVs \\  \hline
$B$ & charging time &$D$& parking time \\  \hline
$V$ & voltage &W& voltage magnitude \\  \hline
$R$ & resistance  &X& reactance\\  \hline
$\blt{p}$ & power allocated to EVs  &$\blt{\Lambda}$& power allocated to nodes \\  \hline
$c_{j}^{\text{max}}$ & maximum charging rate &$M$& node power capacity \\  \hline
\end{tabular}}
\end{table}

\subsection{Network and infrastructure} Consider a radial distribution network described by a rooted tree graph $\mathcal{G}=(\Node,\Edge)$, where $\Node=\{0,1,\ldots, I\}$, denotes its set of nodes (buses) and $ \Edge$ is its set of directed edges, assuming that node $0$ is the root node (feeder). Further, we denote by $\ep_{ik}\in \Edge $ the edge that connects node $i$ to node $k$, assuming that $i$ is closer to the root node than $k$. Let $\Node(k)$ and $\Edge(k)$ be the node and edge set of the subtree rooted in node $k \in \Node$. The active and reactive power consumed by the subtree $(\Node(k),\Edge(k))$ are $P_{\Node(k)}$ and $Q_{\Node(k)}$. The resistance, the reactance, and the active and reactive power losses along edge $\ep_{ik}$ are denoted by $R_{ik}$, $X_{ik}$, $L^P_{ik}$, and $L^Q_{ik}$, respectively. Moreover, $V_i$ is the voltage at node $i$ and $V_0$ is known. At any node, except for the root node, there is a charging station
with $K_i>0$, $i \in \Node \setminus \{0\}$ parking spaces (each having an EV charger). Further, we assume there are $\Type=\{1,\ldots,J\}$ different types of EVs indexed by $j$.

\subsection{Stochastic model for EVs}
Type-$j$ EVs arrive at node $i$ according to a Poisson process with rate $\lambda_{ij}$. If all spaces are occupied, a newly arriving EV does not enter the system, but is assumed to leave immediately.
Each EV has a random charging requirement and a random parking time. These depend on the type of the EV, but are independent between EVs and are denoted by $B_j$ and $D_j$, respectively, for type-$j$ EVs. In queueing terminology, these  quantities are called \textit{service requirements} and \textit{deadlines}, respectively. The joint distribution of $(B_j, D_j)$ is given by a bivariate probability distribution: $F_j(x,y) = P(B_j \leq  x , D_j \leq y)$ for $x,y \geq 0$.
\angb{Our model assumptions are quite standard in the computer-communications setting. For validations in the setting of EVs, we refer to \cite{alizadeh2014, zhang2015}.
The most critical assumption is the assumption that the arrival process is homogeneous, though some of our fluid limit results still hold under time-varying arrival rates. Methods to reduce the case of time-varying arrival rates to fixed arrival rates are explored in \cite{WhittSurvey}.}

Our framework is general enough to distinguish between types. For example, we can classify types according to intervals of ratio of the charging requirement and parking time and/or according to the contract they have. An EV leaves the system after its parking time expires. It may be not fully charged. If an EV finishes its charge, it remains at its parking space without consuming power until its parking time expires. EVs that have finished their charge are called ``\textit{fully charged}''.

\subsection{State descriptor}
We denote by $Q_{ij}(t)\in \{0,1,\ldots, K_i\}$ the total number of type-$j$ EVs at node $i$ at time $t\geq 0$, where $Q_{ij}(0)$ is the initial number of EVs. Thus, $Q_{i}(t):=\sum_{j=1}^{J} Q_{ij}(t)$ denotes the total number of EVs at node $i$. Further, we denote by $Z_{ij}(t)\in\{0,1,\ldots, Q_{ij}(t)\}$ the number of type-$j$ EVs at node $i$ with a not-fully-charged battery at time $t$ and by $Z_{ij}(0)$ the number of vehicles initially at node $i$.
Last, we write $Z_{ij}(\infty)$ or simply $Z_{ij}$ to represent the process in steady-state.

For some fixed time, let $\blt z=(z_{ij}: i\in \Node \setminus \{0\},\ j\in \Type)\in \R_{+}^{I\times J}$ be the vector giving the number of uncharged EVs for all types and nodes. Note that although the vector that gives the number of uncharged EVs should have integer-valued coordinates, we allow real values in order to accommodate fluid analogues later.
Moreover, we assume that EVs receive only active power during their charge and do not absorb reactive power; see  \cite{sojoudi2011optimal} for a justification.

\subsection{Charging control rule}
An important part of our framework is the way we specify how the charging of EVs takes place. Given that the state of uncharged vehicles is equal to $\blt z$, we assume
the existence of a function $\blt p(\blt z)=(p_{ij}(\blt z): i\in \Node \setminus \{0\},\ j\in \Type)$ that specifies the instantaneous rate of power each uncharged vehicle receives.
Moreover, we assume that this function is obtained by optimizing a ``global'' function. Specifically, for a type-$j$ EV at node $i$  we associate a function $u_{ij}(\cdot)$, which is strictly increasing and concave in $\R_{+}$, twice differentiable in $(0,\infty)$ with $\lim_{x \rightarrow 0} u_{ij}'(x)=\infty$. The  charging rate $\blt p(\blt z)$ is then determined by
$\underset{}{\text{max}_{\blt p}}  \sum_{i=1}^{I} \sum_{j=1}^{J} z_{ij} u_{ij}(p_{ij})$
subject to a number of constraints, which take into account physical limits on the charging of the batteries, load limits, and most importantly voltage drop constraints.

Before we describe these constraints in detail, we first provide some comments on this charging protocol. An important example is the choice $u_{ij}(p_{ij}) = w_{ij} \log p_{ij}$, which is known as \textit{weighted proportional fairness}. Note that this scheme assumes the existence of a virtual agent that is capable of optimizing the global function. In practice, this control may be implemented in
a decentralized fashion, using mechanisms such as the \emph{alternating direction } method of multipliers \cite{Chertkov2017, Pinson2017}. It is even possible to come up with decentralized allocation schemes that achieve this control if the functions $u_{ij}(\cdot)$ are unknown, which dictates the use of proportional fairness with a specific choice of the weights $w_{ij}$. For background, we refer to  \cite{kelly1997charging,yi2008stochastic}. A limitation of our formulation is that it does not take into account the remaining time until the deadline
expires and the remaining charging requirement. \angb{Our multiclass framework allows to at least partially overcome this, for example by letting the functions $u_{ij}(\cdot)$ depend on the joint distribution of $(B_j, D_j)$. For example, given a discrete distribution of $(B,D)$, we can classify types $j$ by the ratio $B_j/D_j$, and, in the context of proportional fairness,
modify weights $w_{ij}$ accordingly. Note that it is feasible to communicate an indication of $(B,D)$ by the owner of an EV at the parking lot \cite{arif2016}.}

We now turn to a discussion of the constraints.
 We assume that the highest power that the parking lot at node $i$
 can consume is $M_i>0$ and that the maximum power rate which a type-$j$ EV  can be charged is $c_{j}^{\text{max}}$. That is,
 \begin{equation}\label{eq:nodecon}
    \sum_{j=1}^{J} z_{ij}p_{ij} \leq M_i \quad \text{and}
    \quad
  0 \leq p_{ij} \leq c^{max}_j.
 \end{equation}
We refer to \eqref{eq:nodecon} as load constraints.
\angb{Lines and transformers supply the aggregate load imposed by both homes (and other features) and EV chargers. Loads from homes are usually uncontrolled, and if we would like to explicitly incorporate such loads in our model, we would
need to subtract the load of homes (or an upper bound of that) from the node capacity, and modify the voltage constraints by including uncontrollable loads at each node. To keep the notation manageable, we simply assume here that $M_i$ is the node capacity that can be transferred to EVs. We refer to \cite{ardakanian2013} for more details. }

In addition, we impose voltage drop constraints. These constraints rely on load flow models. Two of these models we consider are described next.

\vspace{-0.3cm}
\subsection{AC voltage model}
We first consider a minor simplification of the full AC power flow equations.
The angle between voltages in distribution networks are small and hence they can be chosen so that the phasors have zero imaginary components
\cite[Chapter 3]{kersting2012distribution}. Under this assumption, Kirchhoff's law \cite[Eq.~1]{low14I} for our model takes the form, for $\ep_{pk}\in \Edge$,
\begin{equation}\label{eq:AC}
  V_p V_k-V_kV_k-P_{\Node(k)} R_{pk}-Q_{\Node(k)} X_{pk}=0,
\end{equation}
where $p\in \Node$ denotes the unique parent of node $k$.
The previous equations are non-linear. Applying the transformation
\begin{align*}
\blt{W}(\ep_{pk})
=\left(
  \begin{array}{cc}
    V_p^2 & V_p V_k \\
    V_k V_p & V_k^2 \\
  \end{array}
\right)
=:
\left(
  \begin{array}{cc}
    W_{pp} & W_{pk} \\
    W_{kp} & W_{kk} \\
  \end{array}
\right)
\end{align*}
leads to the following linear equations (in terms of $\blt{W}(\ep_{pk})$):
\begin{equation}\label{eq:KVL}
W_{pk}-W_{kk}-P_{\Node(k)} R_{pk}-Q_{\Node(k)} X_{pk}=0,\ \ep_{pk}\in \Edge.
\end{equation}
Note that $\blt{W}(\ep_{pk})$ are positive semidefinite matrices (denoted by $\blt{W}(\ep_{pk})  \succeq  0$) of rank one.
The active and reactive power consumed by the subtree $(\Node(k),\Edge(k))$ are given by
\begin{align}
P_{\Node(k)}&= \sum_{l\in \Node(k)}\sum_{j=1}^{J}
z_{lj} p_{lj} +\sum_{l\in \Node(k)}
\sum_{\ep_{ls}\in \Edge(k)}L^P_{ls}, \label{eq:power}\\
Q_{\Node(k)}&=\sum_{l\in \Node(k)}\sum_{\ep_{ls}\in \Edge(k)}L^Q_{ls},\nonumber
\end{align}
where
by \cite[Appendix~B]{carvalho2015critical},
\begin{align*}\label{eq:APL}
L^P_{ls}&=(W_{ll}-2W_{ls}+W_{ss})R_{ls}/(R^2_{ls}+X^2_{ls}),\\
L^Q_{ls}&=(W_{ll}-2W_{ls}+W_{ss})X_{ls}/(R^2_{ls}+X^2_{ls}).
\end{align*}
Note that $W_{kk}$ are dependent on vectors $\blt p$ and $\blt z$. We write $W_{kk}(\blt p,\blt z)$, when we wish to emphasize this. If we use this model to describe voltages, the function
$\blt p(\blt z)$ is then given by
\begin{equation}\label{OP}
\begin{aligned}
& \underset{}{\text{max}_{\blt p, \blt W}}
& & \sum_{i=1}^{I} \sum_{j=1}^{J} z_{ij} u_{ij}(p_{ij}) \\
& \text{subject to}
&&   \eqref{eq:nodecon},\eqref{eq:KVL},
\ \underline{ \upsilon}_i \leq W_{ii} \leq \overline{\upsilon}_i ,\\
&&& \blt{W}(\ep_{ik})  \succeq  0,\ \text{rank}(\blt{W}(\ep_{ik}))=1,\ \ep_{ik}\in \Edge,
\end{aligned}
\end{equation}
for $z_{ij}>0$. If $z_{ij}=0$, then take $p_{ij}=0$. In addition, $0<\underline{ \upsilon}_k \leq W_{00} \leq \overline{\upsilon}_k $ are the voltage limits.
Observe that the optimization problem~\eqref{OP} is non-convex and NP hard due to rank-one constraints. Removing the non-convex constraints yields a convex relaxation, which is a second-order cone program. Further, by Remark~\ref{re:AlCon} (see below) and
 \cite[Theorem 5]{low14}, we obtain that the convex-relaxation problem is exact.

\subsection{Linearized Distflow model}\label{sec:Distflow}
Though the previous voltage model is tractable enough for a convex relaxation to be exact, it is rather complicated. Assuming that
the active and reactive power losses on edges are small relative to the power flows, but now allowing the voltages to be complex numbers, we derive a linear approximation of the previous model, called the \textit{linearized (or simplified) Distflow model} \cite{Baran89}.
In this case, the voltage magnitudes $W_{kk}^{lin}:=|V_{k}^{lin}|^2$ have an analytic expression \cite[Lemma~12]{low14I}:
\begin{equation}\label{eq:Dist}
  W_{kk}^{lin}(\blt p,\blt z)=W_{00}
  -2\sum_{\ep_{ls}\in \mathcal{P}(k)}R_{ls}\sum_{m \in \Node(s)}
  \sum_{j=1}^{J} z_{mj}p_{mj},
\end{equation}
where the $\mathcal{P}(k)$ is the unique path from the feeder to node $k$.
\begin{rem}\label{re:AlCon}
Note that $W_{kk}^{lin}\leq W_{00}$ for all nodes $k$, as we assume that the nodes only consume power, and by \cite[Lemma~12]{low14I} we obtain
$W_{kk}(\blt p,\blt z) \leq W_{kk}^{lin}(\blt p,\blt z)$.
\end{rem}
To derive the representation of the power allocation mechanism $\blt p (\blt z)$ in this setting, one replaces the constraints in~\eqref{OP} by \eqref{eq:nodecon} and
$\underline{ \upsilon}_k \leq W^{lin}_{kk}(\blt p,\blt z) \leq \overline{\upsilon}_k$.

When adding stochastic dynamics, the resulting model is still rather complicated. Even Markovian assumptions yield a multidimensional Markov chain of which the transition rates are governed by solutions of nonlinear programming problems. %
The dynamics of the high-dimensional non-Markovian stochastic process $(Q_{ij}(t), Z_{ij}(t))$, $t\geq 0$, are in general not tractable from a probabilistic viewpoint.
To obtain a Markovian description, we would also have to keep track of all residual parking times and charging requirements, leading to a measure-valued process as in \cite{remerova14}.
Therefore, we consider fluid approximations of $(Q_{ij}(t), Z_{ij}(t))$ in the next section, which are more tractable \angb{and hold under general assumptions.}

\section{Fluid Approximation}\label{sec:FA}
In this section, we develop a fluid approximation for the stochastic model defined in Section~\ref{Sec:model description}, of which the invariant point is characterized through an OPF problem. To do so, we
follow a similar approach as in \cite{remerova14} and \cite{gromoll2006impact}.

The fluid approximation, which is deterministic, can be thought of as
a formal law of large numbers approximation.
More precisely,
consider a family of models as defined in Section~\ref{Sec:model description}, indexed by a scaling parameter $n\in \mathds{N}$. The fluid scaling is given by $\frac{Z^n_{ij}(\cdot)}{n}$. To obtain a non-trivial fluid limit, we choose the following scaling  for the node parameters in the $n^{\text{th}}$ system. The  maximum power that node $i$ can consume is $nM_i$, the arrival rate  is $n\lambda_{ij}$, the number of parking spaces is $nK_i$; all other parameters remain unchanged. A mathematically rigorous justification of this scaling is beyond the scope of this study, and will be pursued elsewhere. If the set of feasible power allocations is polyhedral, the methods from \cite{remerova14} can be applied directly to achieve this justification. This scaling gives rise to the following definition of a fluid model.
\begin{defn}[Fluid model]\label{def:fluid model}
A nonnegative continuous vector-valued function $\blt z(\cdot)$ is a fluid-model solution if it satisfies the functional equations for $i,j\geq 1$
\begin{align*}
  z_{ij}(t)&=z_{ij}(0)
   \Prob\big(B_j^0>\int_{0}^{t}p_{ij}(\blt{z}(u))du,D_j^0\geq t\big)\\
  &+
  \int_{0}^{t}\gamma_{ij}(s) \Prob\big(B_j>\int_{s}^{t}p_{ij}(\blt{z}(u))du,D_j\geq t-s\big)ds,
\end{align*}
where $\gamma_{ij}(t):=\lambda_{ij}\ind {q_i(t)<K_i}$, and $q_i(t) = \sum_{j} q_{ij}(t)$, with
\[
q_{ij}(t) = q_{ij}(0) \Prob \big(D_j^0 \geq t\big)+
  \int_{0}^{t}\gamma_{ij}(s) \Prob\big(D_j\geq t-s\big)ds.
\]
Further, $B_j^0$ and $D_j^0$ are the charging requirement and the parking time for the initial population in the system.
\end{defn}
The time-dependent fluid model solution can be used directly to approximate the evolution of the system at time $t$; e.g., one may take $\E{Z_{ij}(t)}\approx z_{ij}(t)$. This set of equations can be extended to
time-varying arrival rates by replacing $\lambda_{ij}$ by $\lambda_{ij}(t)$.
Also, one can consider schemes in which blocked EVs are not lost, but routed to adjacent parking lots, which leads to further
modifications to $\gamma_{ij}(\cdot)$. The fluid model equations, though still rather complicated, have an intuitive meaning; e.g., the term $\Prob\big(B_j>\int_{s}^{t}p_{ij}(\blt{z}(u))du,D_j\geq t-s\big)$
resembles the fraction of EVs of type-$j$ admitted to the system at time $s$ at node $i$ that are still in the system at time $t$ (for this to happen, their deadline needs to exceed $t-s$ and their service requirement
needs to be bigger than the service allocated, which is $\int_{s}^{t}p_{ij}(\blt{z}(u))du$).
\
We now turn to the behavior of our fluid model in equilibrium, i.e., for $t=\infty$. In this case, we obtain a computationally tractable characterization
through a particular OPF problem. Before we state our main theorem, we introduce some notation. Let
\begin{equation}\label{eq:defGamma}
\gamma_{ij} := \frac{\lambda_{ij}}{\lambda_i} \min \Big\{\lambda_i, K_i \left(\sum_{j'=1}^J \frac{\lambda_{ij'}}{\lambda_i} \E{D_{j'}}  \right)^{-1}\Big\},
\end{equation}
where $\lambda_i:=\sum_{j=1}^{J} \lambda_{ij}$. \angb{Here, $\gamma_{ij}$ can be interpreted as the modified arrival rate of EVs due to finite parking spaces; as not all EVs find an available parking space.}
Furthermore,
define the functions
\begin{equation}\label{eq:defg}
g_{ij}(x):=\gamma_{ij} \E{ \min \{D_j x , B_{j}\}},
\end{equation}
 and the node allocation (the power which type-$j$ EVs consume at node $i$), $\Lambda_{ij}(\blt z):=z_{ij} p_{ij}(\blt z)$. Also, for a random variable $Y$, denote by $\inf(Y)$ the leftmost point of its support.

\begin{thm}[Characterization of invariant point]\label{Th:UnFP}
(i) If  $z^*_{ij}$ is an invariant point for the fluid model, it is given by the solution of the fixed-point equation
\begin{equation}\label{eq:fluid proxy}
z^*_{ij} = \gamma_{ij}
 \E{ \min \{D_j , \frac{ B_{j}} {p_{ij}(\blt z^*)}\}}.
\end{equation}
(ii) Let $\inf{(D_j/B_{j})}\leq 1/c_j^{\text{max}}$.
The solution $\blt z^*$ of \eqref{eq:fluid proxy} is unique and is given by
$z_{ij}^*=\frac{\Lambda_{ij}^*}{g_{ij}^{-1}(\Lambda_{ij}^*)}$, where $\blt \Lambda^*$
is the unique solution of the optimization problem
 \begin{equation}\label{GDOP}
\begin{aligned}
& \underset{}{\text{max}_{\blt{\Lambda},\blt W}}
& & \sum_{i=1}^{I} \sum_{j=1}^{J}  G_{ij} (\Lambda_{ij}) \\
& \text{subject to}
&&W_{ik}-W_{kk}-P_{\Node(k)} R_{ik}-Q_{\Node(k)} X_{ik}=0,\\
&&&\underline{ \upsilon}_i \leq  W_{ii}
\leq \overline{ \upsilon}_i, \blt{W}(\ep_{ik})  \succeq  0,\
\  \Lambda_{ij} \leq M_i,\\
&&&  0 \leq \Lambda_{ij} \leq g_{ij}(c_j^{max}),\ \ep_{ik}\in \Edge.
\end{aligned}
\end{equation}
$G_{ij}(\cdot)$ is a strictly concave function such that $G_{ij}'(\cdot)=u_{ij}'(g^{-1}_{ij}(\cdot))$.
\end{thm}
By \eqref{eq:power}, observe that $W_{kk}$ depends on $\blt{z}$ through the products $z_{ij} p_{ij}(\blt z)$. By the definition of $\blt \Lambda$, we have that $W_{kk}$ depends only on $\blt \Lambda$. That is, the previous optimization problem is indeed independent of the fixed point $\blt z^*$.

The intuition behind (\ref{eq:fluid proxy}) follows from a result in queueing theory known as Little's law \cite{Kleinrock1976},
which says that the expected number of customers equals the arrival rate times the expected sojourn time of a particular customer. The sojourn time of a customer of type-$j$ at node $i$ is the minimum of its deadline and its potential service time. The latter quantity is approximated by the total service time divided by the service rate, which will be approximately constant in equilibrium for large systems
(this is called the snapshot principle in queueing theory \cite{Reiman82}). To arrive at (\ref{GDOP}) the essential idea is to add Little's law (\ref{eq:fluid proxy}) to the set of KKT conditions that characterize
$\blt p ( \blt z)$ and rewrite all equations in such a way that they form the KKT conditions for the problem (\ref{GDOP}). We refer to the Appendix for more details.

\angb{An interpretation of the function to be maximized in Theorem~1 is not straightforward. To provide physical intuition of what is being optimized, we observe that the optimal value of the optimization problem in Theorem 1 can be be rewritten by using the following identity:
\begin{align*}
G_{ij}(\Lambda_{ij}^*)=\gamma_{ij}
\E{ D_j u_{ij}(\min\{p_{ij}^*,\frac{B_j}{D_j}\})},
\end{align*}
with $p_{ij}^*= \Lambda_{ij}^*/z_{ij}^*$. For the calculations leading to this identity, we refer to
Appendix~\ref{Interpretation of the objective function}. Note that $G_{ij}(\cdot)$ depends naturally on the joint distribution of $B$ and $D$. Further,
the ratio $\frac{B_j}{D_j}$ can be interpreted as the desired power of type-$j$ EVs and $p_{ij}^*$ is the actual power it receives in the long term. So, we pick $\blt{p}^*$ that maximizes the long term utility
$\sum_{i=1}^{I} \sum_{j=1}^{J}
\gamma_{ij}
\E{ D_j u_{ij}(\min\{p_{ij},\frac{B_j}{D_j}\})}$.
This interpretation heuristically follows by using a version of Little's law: EVs of type $j$ arrive at node $i$ at rate $\gamma_{ij}$, stay for time $D_j$, and during their stay they are interested in a rate not exceeding $B_j/D_j$, so the utility of rate $p_{ij}$ is in equilibrium
$u_{ij}(\min\{p_{ij},\frac{B_j}{D_j}\})$.}

To show that a fluid model solution converges towards the equilibrium point as $t\rightarrow\infty$, we also assume that the power rate that an EV gets does not increase if a newly arriving EV enters the system.
\begin{defn}
An allocation mechanism is called ``monotone'' if\quad $0<\blt{y}\leq \blt{z}$ implies that
$p_{ij}(\blt{y})\geq p_{ij}(\blt{z}) $.
\end{defn}
By adapting techniques from \cite{BorstEgorovaZwart2014}, it can be shown that this assumption holds in the linearized Distflow model.
%
\begin{prop}\label{eq:SSL}
If the allocation mechanism is monotone and $K_i=\infty$, we have that $z_{ij}(t) \rightarrow z^*_{ij}$, as $t \rightarrow \infty$.
\end{prop}
\angb{Using the continuous mapping theorem \cite{chen2001fundamentals} and extending results from \cite{kang2015}, \cite{gromoll2009}, and  \cite{remerova14}, we conjecture this result to be true without imposing the condition $K_i=\infty$ and also allowing for general distributed arrival process, but the mathematical techniques required to establish this are beyond the scope of this paper.}

\angb{As we already mentioned the original stochastic system is not tractable even for computational purposes. In Appendix~\ref{Additional numerical examples},
we compare the fluid approximation with its original stochastic model and also compare across load flow models. Our computational study, though certainly not exhaustive, indicates that our fluid models can accurately predict the performance of a specific control rule.}

\section{\angb{Case study}}\label{sec:case Study}
In this section, we apply our results to the SCE 47-bus network; \cite{low14}. Node 1 is the feeder and nodes 13, 17, 19, 23 and 24 (with white color in Figures~\ref{fig:PerCl1}, \ref{fig:PerCl2}, \ref{fig:PerSNM}, and \ref{fig:SPM}) are photovoltaic generators, and we removed them from the network. We present results for both the aggregated system and the individual nodes.
First, we study a non-Markovian system with multiple type of EVs where the charging and parking times are dependent. Then, we move to a single class non-Markovian model. Last, we present results for a Markovian model where we also compare the SCE 47-bus with the SCE 56-bus network. In all cases, we consider the proportional fairness allocation, i.e., $u_{ij}(x)=\log(x)$. Further, we fix $\blt{M}=8$, $c^{\text{max}}=1$, and allow the voltage drop to diverge at most 10\% from its nominal. Last, we use MATLAB and CVX to solve the semidefinite optimization problems. It takes approximately 3 minutes for the SCE 47-bus network and 5 minutes for the SCE 56-bus network.

\subsection{Non-Markovian model with multiple type of EVs}
\label{sec:nonMmulti}
The first case we study is a non-Markovian network with multiple type of EVs with dependent charging and parking times. In particular, we consider the case where
parking and charging times satisfy $B_j/D_j=\theta_j$, where $\theta_j>0$ is deterministic. Here, $\theta_j$ can be interpreted as the desired power of type-$j$ EV. In this case, we have that
\begin{align*}
G_{ij}(x)=\gamma_{ij} \E{D_j}\log(\min\{x,\theta_j\}).
\end{align*}

\begin{figure}[!h]
		\includegraphics[width=.8\linewidth]
		{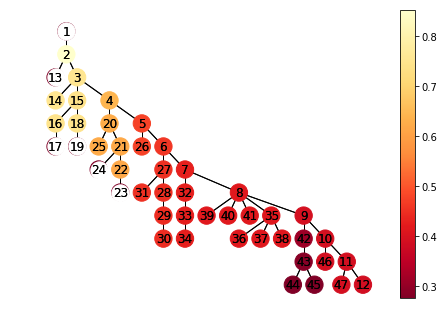}
	\caption{Percentage of power type-$1$ EVs receive}
\label{fig:PerCl1}
\end{figure}

\begin{figure}[!h]
		\includegraphics[width=.8\linewidth]
		{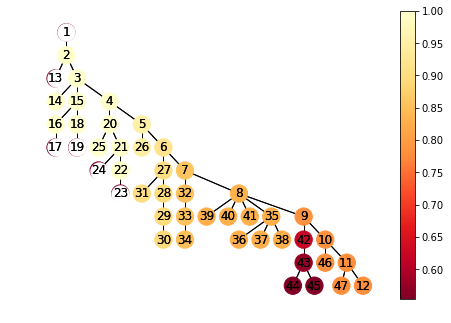}
	\caption{Percentage of power type-$2$ EVs receive}
\label{fig:PerCl2}
\end{figure}

Consider an overloaded system taking
$\blt{K}=1$, $\lambda=1.2$, $\E{D_j}=1$ and two type of EVs with
arrival rate $\lambda_{i1}=0.4 \lambda$ and $\lambda_{i2}=0.6 \lambda$, respectively. In Figures~\ref{fig:PerCl1} and \ref{fig:PerCl2}, we see the percentage of the desired power received (i.e., $p_{ij}/\theta_j$)  by both types of EVs and all the nodes for $\theta_1=0.02$ and $\theta_2=0.01$. As expected, the power grid is able to provide better quality of service to EVs with lower $\theta_j$. Type-$2$ EVs receive power very close to their preference (i.e., $\theta_2$) at nodes close to the feeder and it reduces for nodes away from the feeder. On the other hand, type-$1$ EVs receive less power than their preference and it can significantly drop for nodes away from the feeder. Note that it is not feasible for the power grid to allocate power $\theta_j$ to all the EVs as this would cause a significant voltage drop.

\subsection{Non-Markovian model with single type of EVs}
\begin{figure}[!h]
		\includegraphics[width=.8\linewidth]
		{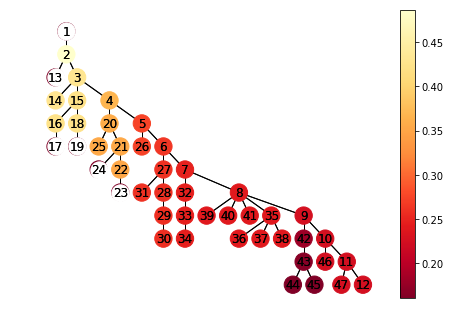}
	\caption{Ratio $\frac{\Lambda}{\E{\Theta}}$ in case of $B=\Theta D$}
\label{fig:PerSNM}
\end{figure}

Here, we study a relation between charging and parking times in the form
$B=\Theta D$, where $\Theta$ is a discrete random variable with $\Prob(\Theta=\theta_i)=q_i$, for $\theta_i>0$. For simplicity, we take that $\Theta\in \{\theta_1,\theta_2\}$ and without loss of generality we take $\theta_1< \theta_2$. In this case,  we have that
\begin{align*}
  G_i(x)=\gamma_i \E{D}(\log(\min\{x,\theta_1\})
+q_2\log(\min\{x,\theta_2\}\\
-\gamma_i \E{D}\theta_1q_1).
\end{align*}
In Figure~\ref{fig:PerSNM}, we see the ratio $\frac{\Lambda}{\E{\Theta}}$ for an overloaded system, with $\blt{K}=1$, $\E{D}=1$, and $\lambda_i=1.2$. Moreover, we choose $\theta_1=0.001$, $\theta_2=0.02$, $q_1=0.1$, and $q_2=0.9$. The aggregated success probability (i.e., the probability an EV leaves the network with fully charged battery) in this case is approximately 0.1. Every individual node receives less power (almost the half) compared to the case studied in Section~\ref{sec:nonMmulti}. This can be explained as follows. In the second case, all EVs desire power $\theta_2$ with high probability compared to the case studied in Section~\ref{sec:nonMmulti} where only type-$2$ EVs have this preference.

\subsection{Markovian model}
\begin{figure}[!h]
	\centering
		\includegraphics[width=.8\linewidth]
		{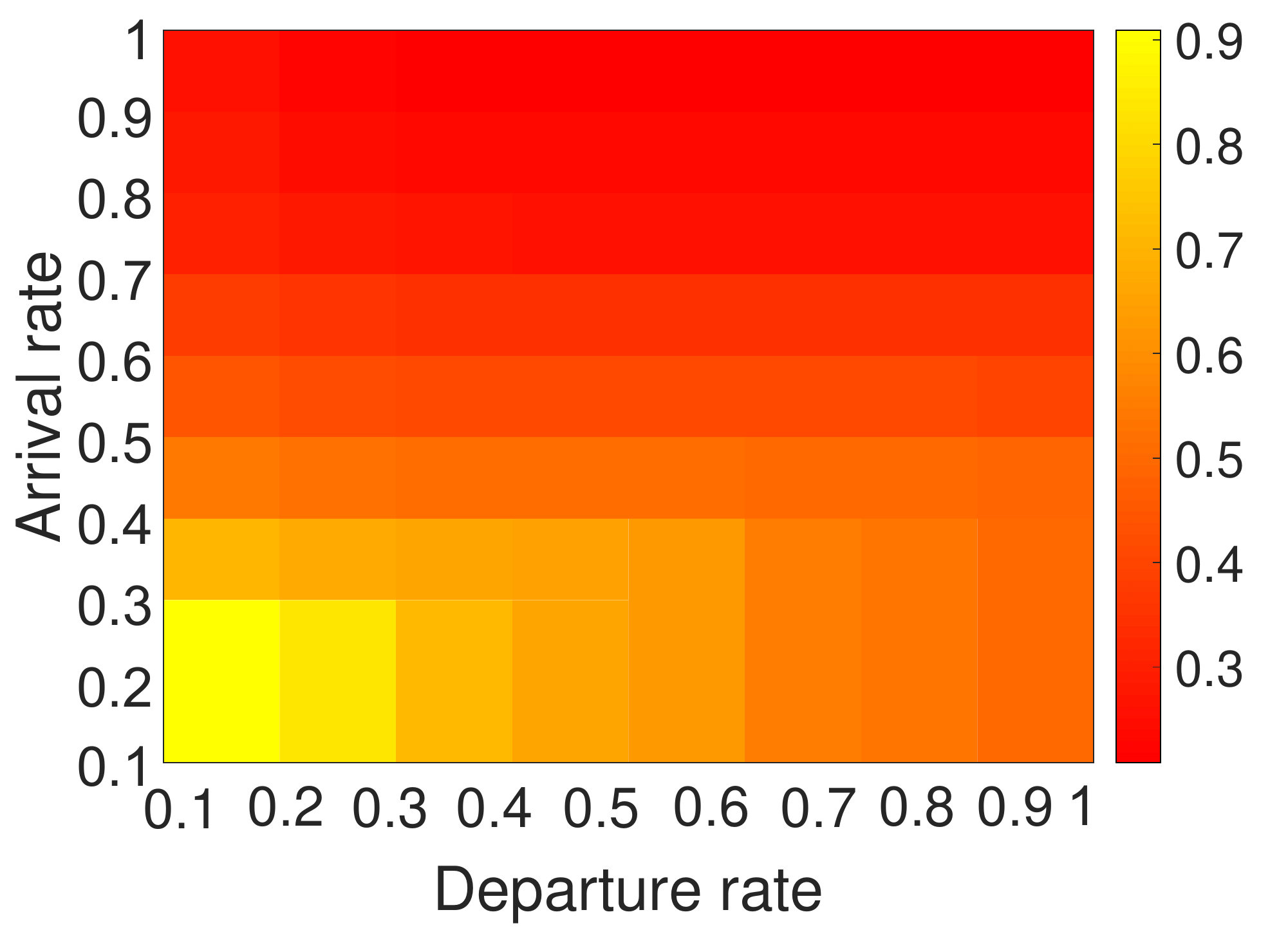}
\caption{Aggregated success probability in SCE 47-bus network}
\label{fig:ASPM}
\end{figure}

\begin{figure}[!h]
	\centering
		\includegraphics[width=.8\linewidth]
		{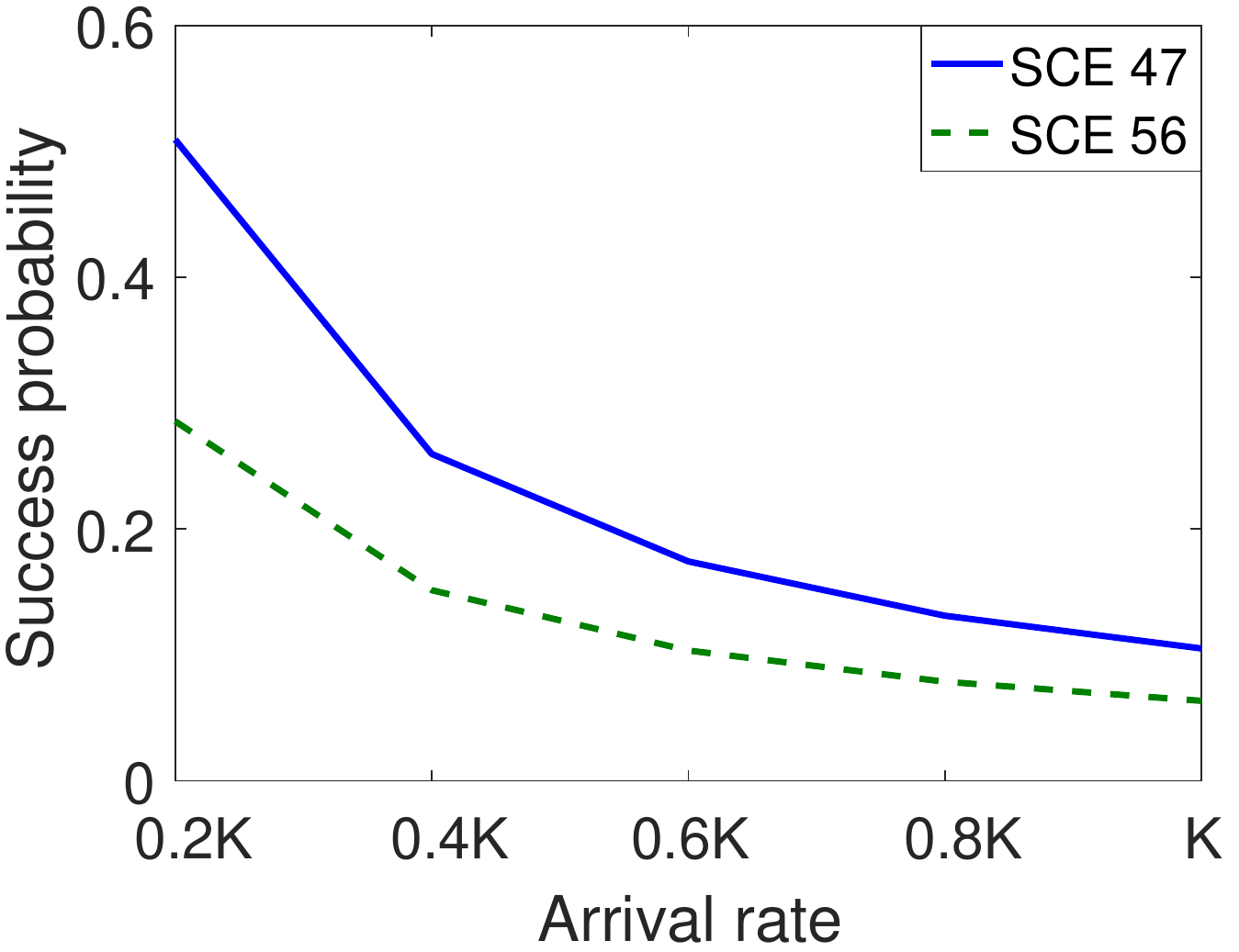}
	\caption{Aggregated success probability in SCE 47-bus and SCE 56-bus networks}
\label{fig:ASCP4756}
\end{figure}

Here, we assume that the charging time $B$ and parking time $D$ of the EVs are independent exponential distributed. EVs arrive in the system with arrival rate $\lambda$ and choose a node uniformly. We first focus on the whole system fixing  $\blt{K}=1$ and $\E{B}=1$. The aggregated success probability in the SCE 47-bus  network for different values of arrival rate $\lambda$ and departure rate $1/\E{D}$ is shown in Figure~\ref{fig:ASPM}.
We observe that the aggregated success probability decreases as more EVs arrive in the system and as the EVs become more impatient (i.e., as the departure rate increases and as a result the parking times become larger). We comment that for small arrival rate and small departure rate the power grid is able to charge the EVs at the maximum rate and since $c^{\text{max}}=1$, the success probability is given by $\Prob{(D>B)}=\frac{1/\E{D}}{1/\E{D}+1/\E{B}}$. Hence, the previous formula gives an upper bound of the success probability for any $\E{B}$ and $\E{D}$. In Figure~\ref{fig:ASCP4756}, we compare the aggregated success probability as a function of the arrival rate between SCE 47-bus and SCE 56-bus networks. Here, we choose $\blt{K}=10$, $\E{B}=0.2$, and $\E{D}=1$. The behavior of the SCE 47-bus and SCE 56-bus networks i similar. However, the success probability in the SCE 56-bus network is smaller as it contains more buses. Moreover, in both networks the probability an EV leaves the system with a full battery is decreasing as a function of the arrival rate.

\begin{figure}[h!]
	\centering
		\includegraphics[width=.8\linewidth]
		{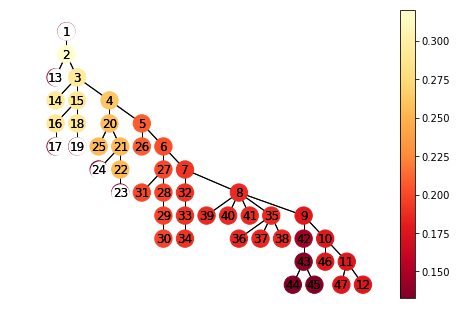}
	\caption{Individual success probabilities}
\label{fig:SPM}
\end{figure}

Next, we move to the behaviour of the individual nodes. In Figure~\ref{fig:SPM}, we can see the success probability for all the nodes for $\lambda=1$ and $\E{B}$=$\E{D}=1$. Not surprisingly the success probability decreases as we go away from the feeder as the transferred power is reduced.

\section{Explicit results under linearized Distflow model and a line topology}
\label{Additional examples under a line topology}
\angb{
Although Theorem~\ref{Th:UnFP} provides a computationally attractive procedure to analyze the long-term performance of a class of EV charging schemes in a general setting, it does not provide mathematically explicit formulas. The goal of the present section is to present a number of explicit examples, at the expense of making additional assumptions. Specifically, we assume the case of a line network, and will only focus on the linearized DistFlow model. In addition, we will only consider the voltage drop constraints.}

\angb{
We investigate four cases, one special case focusing in the original stochastic model, and the other three cases focusing on the fluid model.}

\subsection{An allocation mechanism with the product-form property}
\label{Sec3}
Despite the complexity of our stochastic model, we are able to identify a special case for which the entire network behaves like a multiclass
processor-sharing queue, of which the invariant distribution
is explicit, and for which even time-dependent properties are known \cite{FralixZwart11, kelly2014stochastic}.

We take $J=1$ for convenience and drop all indices $j$ from the notation in this section. For every node $i$, let
$\SR_i=\sum_{\ep_{ls}\in \mathcal{P}(i)}R_{ls}$,
 $\rho_i = \lambda_i \E{B} \SR_i/\delta$ and
$\rho=\sum_{i=1}^{I}\rho_i$, where
$\delta=\frac{W_{00}-\underline{\upsilon}}{2}$ and $\mathcal{P}(i)$ denotes the unique path from node $i$ to feeder.

\begin{prop}\label{prop:product form}
Assume $K_i=\infty$ for all $i$. If the power allocation rule is proportional fairness, i.e.,
$u_i(p_i)=\log p_i$, then, for every $\blt n\in \mathds{N}^I_+$,
\begin{equation}\label{eq:PS}
\lim_{t\rightarrow\infty} \Prob\big(\blt{Z}(t)=\blt{n}\big)
= (1-\rho)(\sum_{i=1}^{I}n_i)!\prod_{i=1}^I \frac{\rho_i^{n_i}}{n_i!},
\end{equation}
provided $\rho<1$.
\end{prop}
Note that this result is valid for arbitrary distributions of the charging requirements and as such it provides an insensitivity property for the original stochastic model.
This result is another exhibition of the appealing nature of proportional fairness, which has also been shown to give similar nice properties in communication network models
\cite{BonaldProutiere2002, Bonald2006, borst05,  kelly2014stochastic}.
\angb{What is novel in this setting is that the service requirements are stretched with a factor $\SR_i$.}

The proof of this theorem follows from a similar argument as in \cite{BonaldProutiere2002}, making a connection with the class of Whittle networks, by showing that a specific local balance property called \textit{balanced fairness} is satisfied. We explain this procedure for the case of exponential charging times in Appendix~\ref{proofs}.

\subsection{An explicit solution of the time dependent fluid model equations}
In general, the time-dependent fluid model equations (see Definition~\ref{def:fluid model}) need to be solved numerically, for example by using Picard iteration.
In some cases, the set of equations can also be solved explicitly, as the next proposition illustrates.
\begin{prop}\label{prop:explicitZ}
 Assume that $K_i=\infty$, $J=1$, and the charging requirements and parking times are exponential and independent. If $u(p_i)=\log(p_i)$, we have that
 \begin{equation*}\label{eq:Explicit}
z_i(t)=z_i^*+(z_i(0)-z_i^*)e^{-t/\E{D}},
\end{equation*}
with
\begin{equation*}
\label{eq:FA0}
z_i^*=\E{D} \left(\lambda_i-\Lambda_i^*/\E{B}\right), \hspace{0.2cm} \Lambda_i^*= \frac{\lambda_i\delta}{\sum_{k=1}^{I}\SR_k\lambda_k}.
\end{equation*}
\end{prop}
The previous proposition continues to hold for $K_i$ big enough such that the parking lots are never full.

\subsection{A fairness property}
Here, we identify a special case of our setting in which $\blt p(\blt z)$ is explicit and leads to a fair allocation of power to all users, in the sense
that the charging rate of a battery does not depend on the location where it is parked.

An important question is under what assumptions all EVs in the system are charged at the same power rate. The following proposition gives a partial positive answer to this question. Let
$\SR_i=\sum_{\ep_{ls}\in \mathcal{P}(i)}R_{ls}$ and $\delta=\frac{W_{00}-\underline{\upsilon}}{2}$.
\begin{prop}[Fairness property]\label{prop:load balancing}
Let $u_{ij}(p_{ij})=w_{ij} u(p_{ij})$.
If $w_{ij}=\SR_i,$ then we have that $p_{ij}(\blt z)=p(\blt z)>0$. Moreover,
$p(\blt z)= \delta
\Big(\sum_{i=1}^{I}\SR_i\sum_{j=1}^{J}z_{ij}\Big)^{-1}$.
\end{prop}
A similar result can be shown for general trees, under the assumption that the root node has only one child. This choice of weights seems to lead to a low over-all fraction of fully charged cars: in Section~\ref{Designing control schemes for a line}, we show that that fraction is optimized by having the weights $w_i$ decrease with $i$.

\subsection{Designing control schemes for a line}
\label{Designing control schemes for a line}
A natural question is how to choose suitable functions $u_{ij}(\cdot)$ in order to optimize the behavior of the system. In this section, we illustrate the scope of this question
by focusing
the weighted proportional fairness allocation and single type of EVs.
The main goal is to describe a procedure of choosing
the weights $w_i$, $i\geq 1$, such that the log-run fraction of EVs that gets successfully charged is  maximized.  This is a particularly relevant problem in an overloaded regime and to this end we make the assumption $\sum_{i=1}^{I} \E{B}\gamma_i \SR_i > \delta$. 

The solution of ~\eqref{GDOP} is given by
$\Lambda_i^* = g_i\left(\frac {w_i}{h \SR_i }\right)$,
where $g_i(\cdot)$ is given by \eqref{eq:defg}. Under the assumptions of this section, \eqref{GDOP} has only a single constraint. Its associated Lagrange multiplier $h>0$ satisfies
$\delta = \sum_i \SR_i \Lambda_i^* = \sum_i \SR_i g_i(w_i/(h \SR_i ))$.
It follows by Theorem~\ref{Th:UnFP} that
$z_i^*=\frac{h \SR_i\Lambda_i^*}{w_i}$.
The probability an EV gets successfully charged at node $i$ is given by
$\Prob(D>B z_i^*/\Lambda_i^*) = \Prob(w_i D>h \SR_i B )$. Note that the last probability depends on the Langrange multiplier which depends on the weights in an intricate fashion. To overcome this difficulty we scale $\blt w$ such that $h(\blt w)=1$. This can always be achieved as multiplying all weights with the same factor scales $h$ accordingly but leads to the same solution $\blt{\Lambda}^*$. Thus, we can  formulate the problem of optimally selecting weights as follows:
 \begin{equation}\label{OP:OptimalW}
\begin{aligned}
& \underset{}{\text{max}_{\blt{w}}}
& & \sum_{i=1}^{I} \gamma_i \Prob(w_i D> \SR_i B )  \\
& \text{subject to}
&& \sum_{i=1}^{I} \gamma_i \E{\min\{w_iD,\SR_i B\} }\leq \delta.
\end{aligned}
\end{equation}
This problem can be transformed into a resource allocation problem as considered in \cite{katoh98}. It is in general non-convex and its solution depends on the joint distribution of $B$ and $D$. In the sequel, we consider the special case $B=H D$, with $H$ a random variable independent of $D$. Below, we study some special cases of $H$. In turns out that the case of $B/D$ equals to a deterministic is relevant in the setting of distributionally robust optimization.

\subsubsection{$H$ is random  with decreasing hazard rate}
Here, we assume that $H$ is a random variable where its distribution has decreasing hazard rate. The latter ensures that transformation of \eqref{OP:OptimalW} leads to a convex optimization problem. A well-known example of such a distribution is
a Pareto distribution with parameters $ a>1$ and $\kappa>0$, i.e., $\Prob(H>x)=(\frac{\kappa}{x+\kappa})^a$, $x\geq0$. Note that $\E{H}= \kappa/(a-1)$.
Setting $c_i=w_i/ \SR_i$, \eqref{OP:OptimalW} reduces to
 \begin{equation*}
\begin{aligned}
& \underset{}{\text{max}_{\blt{c}}}
& & \sum_{i=1}^{I}\gamma_i [1-\big(\frac{\kappa}{c_i+\kappa}\big)^a]  \\
& \text{subject to}
&& \sum_{i=1}^{I} \frac{ \E{D}\kappa\gamma_i \SR_i}{1-a} [\big(\frac{\kappa}{c_i+\kappa}\big)^{a-1}-1] \leq \delta.
\end{aligned}
\end{equation*}
After setting $y_i=(\frac{\kappa}{c_i+\kappa})^{a-1}$, we obtain $c_i=\frac{\kappa}{y_i^{1/(a-1)}}-\kappa$ and
 \begin{equation*}\label{OP:pareto}
\begin{aligned}
& \underset{}{\text{max}_{\blt{y}}}
& & \sum_{i=1}^{I}\gamma_i (1-y_i^{a/(a-1)}) \\
& \text{subject to}
&& \sum_{i=1}^{I}\frac{ \E{D} \kappa\gamma_i \SR_i}{1-a} (y_i-1) \leq \delta,
&& 0 \leq y_i\leq 1,
\end{aligned}
\end{equation*}
which is a convex problem.

\subsubsection{$H$ is deterministic}
Let $H$ be deterministic and strictly positive.
In this case, setting $c_i=w_i/ \SR_i$, \eqref{OP:OptimalW} reduces to
 \begin{equation*}
\begin{aligned}
& \underset{}{\text{max}_{\blt{c}}}
& & \sum_{i=1}^{I}\gamma_i\ind{H\leq c_i}  \\
& \text{subject to}
&& \sum_{i=1}^{I} \E{D}\gamma_i \SR_i\min\{c_i,H\}\leq \delta.
\end{aligned}
\end{equation*}
There exists a non-empty set $\Node^* \subseteq \Node \setminus \{0\}$ such that
$c_i<H$ for each $i\in \Node^*$. To see this observe that if $\Node^*=\emptyset$, then
\begin{align*}
\delta \geq \sum_{i=1}^{I} \E{D}\gamma_i \SR_i\min\{c_i,H\}
=\sum_{i=1}^{I} \E{B}\gamma_i \SR_i,
\end{align*}
which contradicts the overload assumption. Moreover, observe that the nodes in $\Node^*$ do not increase the objective function and so we can set $c_i=0$ for $i\in \Node^*$. On the other hand, if $i\notin \Node^*$, then $c_i\in[H,\infty)$ and minimization of the constraint leads to $c_i=H$. So, the last optimization problem is equivalent to
 \begin{equation}\label{OP:KP}
\begin{aligned}
& \underset{}{\text{max}_{\blt{x}}}
& & \sum_{i=1}^{I}\gamma_i x_i  \\
& \text{subject to}
&& \sum_{i=1}^{I} \E{D}\gamma_i \SR_iH x_i\leq \delta,
&& x_i\in \{0,1\},
\end{aligned}
\end{equation}
which is the well-known  knapsack problem.

We now show that the case of deterministic $H$ is the worst-case situation in overload by showing that the value of~\eqref{OP:KP} is a lower bound for the value of~\eqref{OP:OptimalW}. Assume w.l.o.g.\  $\E{H}=1$ and
use the Markov inequality to obtain
$\sum_{i=1}^{I} \gamma_i
\Prob(H<c_i) \geq \sum_{i=1}^{I} \gamma_i (1-1/c_i) \ind{c_i>1} $. In addition, an upper bound for the constraint is given by
$\sum_{i=1}^{I}\E{D} \gamma_i \SR_i \E{ \min\{H, c_i\}} \leq
\sum_{i=1}^{I}\E{D} \gamma_i \SR_i \ind{c_i>0}$
 due to the fact that
 $\E{\min\{H, c_i\}} = 0$ if $i\in \Node^*$ and
 $\E {\min\{H,c_i\} }\leq \E{H}=1$, otherwise.
These bounds leads to the problem 
\begin{equation*}\label{OP:LB}
\begin{aligned}
& \underset{}{\text{max}_{\blt{c}}}
& & \sum_{i=1}^{I} \gamma_i (1-1/c_i) \ind{c_i>1}   \\
& \text{subject to}
&& \sum_{i=1}^{I}\E{D} \gamma_i \SR_i \ind{c_i>0}\leq \delta.
\end{aligned}
\end{equation*}
It is clear that the feasible set of the last problem is a subset of the feasible set of~\eqref{OP:OptimalW}. Taking now $c_i=\infty$ if $c_i>0$ yields \eqref{OP:KP} with $H=1$.

The numerical results provided in next section confirm that the case of deterministic $H$ indeed leads to a lower fraction of successful charges than the case
of a random $H$.  This in itself results in a selection of weights where nodes far away from the feeder will not receive any power if $H$ is deterministic, while this is not the case
when $H$ is random.

\subsubsection{Numerical validations for optimal weights}\label{Optimal weights}
Now, we move to the choice of optimal weights. We take ten nodes and parameters of the system such that
$\sum_{i=1}^{10} \frac{\E{B}\gamma_i \SR_i }{\delta}=1.17$; i.e., the system is $17\%$ overloaded. We present three examples.
First, for $H$ deterministic and equal to 1. Second, for $H$ being Pareto distributed with $\E{H}=1$ and two shape parameter $a$, namely $a=1.1$ and $a=3$.

For the deterministic case, we have that the optimal weights are
(0.0610, 0.0915, 0.1110, 0.1268,  0.1390, 0.1488, 0.1573, 0.1646, 0, 0), and the number of successful charges per unit of time equals 8.

For the Pareto case, the optimal weights $w_i$, for $a=1.1$ are given by the vector
 (0.1105, 0.1031, 0.1006, 0.0994, 0.0987, 0.0981, 0.0978, 0.0975, 0.0972, 0.0971), and for $a=3$,
 (0.1293, 0.1173, 0.1093, 0.1033, 0.0985, 0.0945, 0.0911, 0.0881, 0.0854, 0.0830).
The number of successful charges per unit of time equals $10$ and $9.49$, respectively.

It turns out that the Pareto case works better and also the more variability (smaller $a$), the better the system performance. For $H$ being Pareto distributed, the weights are decreasing and all the nodes receive strictly positive power in contrast to the deterministic case. In the latter case, the weights are increasing and some nodes do not receive power. This can be explained intuitively as follows:
the nodes away from the feeder receive more power in order the successful charges in the system being maximized.

\section{Concluding remarks}
This paper has proposed a queueing network model for electric vehicle charging. The main result is a fluid approximation of the number of uncharged vehicles in the system, which is derived by combining key ideas from network utility
maximization (using a class of utility-maximizing scheduling disciplines),  from queueing theory (Little's law and the snapshot principle), and from power systems engineering (AC and linearized Distflow load flow models).
Our fluid approximation explicitly captures the interaction between these three elements, as well as physical network parameters and can be computed using convex programming techniques.
Our approach can easily be extended to impose other reliability constraints, such as line limits, and our fluid approach can be extended to deal with superchargers, by allowing additional queues at each station.

We focused on the specific class of weighted proportional fairness protocols. Our optimization framework allows for a further comparison between charging protocols, which is a natural next step for further research.
Another important problem is to  extend our model to allow for batteries to discharge and to include other features, such as
smart appliances/buildings/meters, rooftop solar panels, and other sources of electricity demand.  This naturally leads to various questions about pricing schemes.
We think that our characterization of the performance of the system for a fixed control in terms of the  OPF problem (\ref{GDOP}) can be extended in such directions (given the universal applicability of Little's law) and can be
a useful starting point to design economic mechanisms.

A full mathematical examination of our system is beyond the scope of the present work. In a follow-up work, we will rigorously show that $z_{ij}^*$ is a good approximation of $Z_{ij}(\cdot)$ by developing
a fluid limit theorem using the framework of measure-valued processes, extending the framework of \cite{remerova14}.

\bibliographystyle{IEEEtran}
{\footnotesize	
\bibliography{../Mybibliography_Angelos_Ab}
}

\appendices
\section{\angb{Interpretation of the objective function}}
\label{Interpretation of the objective function}
Recall that
$G_{ij}'(\cdot)=u_{ij}'(g^{-1}_{ij}(\cdot))$, which can be written as
\begin{equation}\label{SCP}
G_{ij}(\Lambda_{ij}^*)= \int_{y_0}^{\Lambda_{ij}^*}
u_{ij}'(g^{-1}_{ij}(y))dy.
\end{equation}
Choose $y_0$ such that $u_{ij}(g^{-1}_{ij}(y_0))=0$. Defining $x=g^{-1}_{ij}(y)$ and recalling $p_{ij}^*=g^{-1}_{ij}(\Lambda_{ij}^*)$, \eqref{SCP} takes the form
\begin{align}\label{eq:G}
G_{ij}(\Lambda_{ij}^*)= \int_{g^{-1}_{ij}(y_0)}^{g^{-1}_{ij}(\Lambda_{ij}^*)}
u_{ij}'(x)dg_{ij}(x)
=&
 \int_{g^{-1}_{ij}(y_0)}^{p_{ij}^*}
u_{ij}'(x)dg_{ij}(x)\\ \nonumber
 =&
 \int_{g^{-1}_{ij}(y_0)}^{p_{ij}^*}
g_{ij}'(x)du_{ij}(x).
\end{align}
Now, note that by \eqref{eq:defg}, we have that
\begin{align*}
g_{ij}'(x)=& \lim_{\epsilon \rightarrow 0}
 \frac{\gamma_{ij}}{\epsilon}
 \E{\min\{D_j (x+\epsilon) , B_{j}\}-\min\{D_j x, B_{j}\}}\\
 =& \gamma_{ij} \E{ D_j \ind{D_jx<B_j} }.
\end{align*}
Applying the last equation in \eqref{eq:G}, we have that
\begin{align*}
G_{ij}(\Lambda_{ij}^*)=&
 \gamma_{ij} \int_{g^{-1}_{ij}(y_0)}^{p_{ij}^*}
\E{ D_j \ind{D_jx<B_j} } du_{ij}(x)\\
=&\gamma_{ij}
\E{ D_j \int_{g^{-1}_{ij}(y_0)}^{\min\{p_{ij}^*,\frac{B_j}{D_j} \} } du_{ij}(x)}\\
=&\gamma_{ij}
\E{ D_j u_{ij}(\min\{p_{ij}^*,\frac{B_j}{D_j})}.
\end{align*}
\section{Additional numerical examples}
\label{Additional numerical examples}
In this section, we give an additional numerical illustration of our results.

\subsection{Stochastic model and fluid approximation}
First, we examine the stochastic model and the fluid approximation.
We consider the case $J=1$, and assume $I=2$. 
Fix $V_0=1$ in the per-unit system, take $R_{01}=X_{01}=0.01$, $R_{12}=X_{12}=0.005$, and $\underline\upsilon_1=\underline\upsilon_2=(0.9)^2$. We consider a weighted proportional fairness allocation with weights
$w_i =\SR_i=\sum_{\ep_{ls}\in \mathcal{P}(i)}R_{ls}$.
Further, the charging requirements and the parking times are exponentially distributed with unit mean, and we take $K_1=K_2=K$ in all our computations, and only take voltage constraints into consideration.

As a first experiment, we consider this model with only voltage drop constraints.
We take $\lambda_i=1.2K_i$ and compare $\E{Z_i}$ between the AC and Distflow allocation mechanisms, as well as their fluid approximations $z_i^*$.
We compare the two power flow models for the stochastic model computing the expected number of uncharged EVs at any node; see Table~\ref{table1}.
Note that the expectations for the AC model are higher as this model takes into account power losses leading to lower service rates.

\begin{table}[!h]
\centering
\caption{}
\label{table1}
\begin{tabular}{|c| |c| c| c| }
 \hline
 $K$& $\E{Z_1},\E{Z_2}$ (Distflow)  & $\E{Z_1},\E{Z_2}$ (AC)  &Rel. error \\
 \hline
 $10$ &  $4.5336, 4.6179  $ &4.6801, 4.6924&3.13 \%, 1.58 \%\\
 $20$ &  $14.0174, 14.0385$ &14.1725, 14.1948& 1.09 \%, 1.10 \%\\
 \hline
\end{tabular}
\end{table}
Next, we evaluate the fluid approximation for the two load models; see Table~\ref{table2}.
Observe again that the number of uncharged EVs for the AC model is higher.
The relative error between the two load flow models is similar to what we saw in Table~\ref{table1}.
\begin{table}[!h]
\centering
\caption{}
\label{table2}
\begin{tabular}{|c| |c| c| c| }
 \hline
 $K$& $z_1^*,z_2^*$ (Distflow)  & $z_1^*,z_2^*$  (AC) &Rel. error \\
 \hline
 $10$ &  $4.5769,4.5769  $ &$4.7356,4.7513$& $3.35 \%, 3.67 \%$ \\
 $20$ &  $14.0300, 14.0300$ &$14.1849, 14.2069$& $1.09 \%, 1.25 \%$ \\
 $30$ &  $23.6820,23.6820$ &$23.8357, 23.8597$& $0.64 \%, 0.74 \%$ \\
 $40$ &  $33.4293, 33.4293$ &$33.5823,33.6073$& $0.45 \%, 0.53 \%$ \\
 $50$ &  $43.2330 ,43.2330$ &$43.3857,43.4112$& $0.35 \%, 0.41 \%$ \\
 \hline
\end{tabular}
\end{table}
Not surprisingly,
The original stochastic model is not numerically tractable for high values of $K$. For low values of $K$, the difference between the stochastic and fluid models is small:
\begin{table}[!h]
\centering
\caption{}
\label{table3}
\begin{tabular}{ |p{0.5cm}||p{2cm}|p{2cm}|  }
\hline
 \multicolumn{3}{|c|}{Relative error of fluid approximations } \\
 \hline
 $K$& Distflow &AC \\
 \hline
 $10$ &  $0.95 \%, 0.86 \% $ &$1.18 \%, 1.25 \% $\\
 $20$ &  $0.09 \%, 0.06 \% $ & $0.08 \% , 0.08 \% $\\
 \hline
\end{tabular}
\end{table}

\begin{figure}[!h]
	\centering
		\includegraphics[width=.45\linewidth]
		{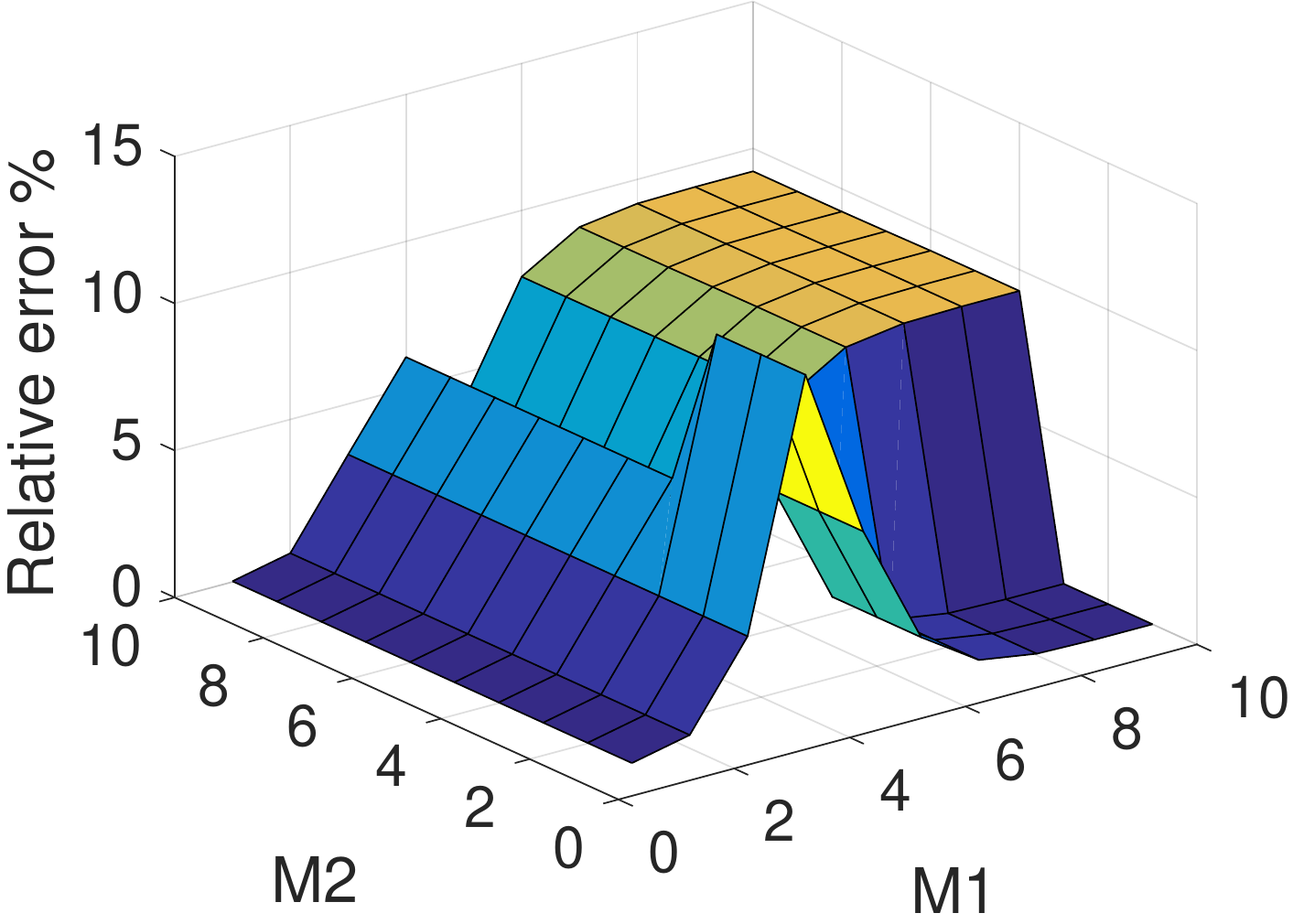}
		\includegraphics[width=.45\linewidth]
		{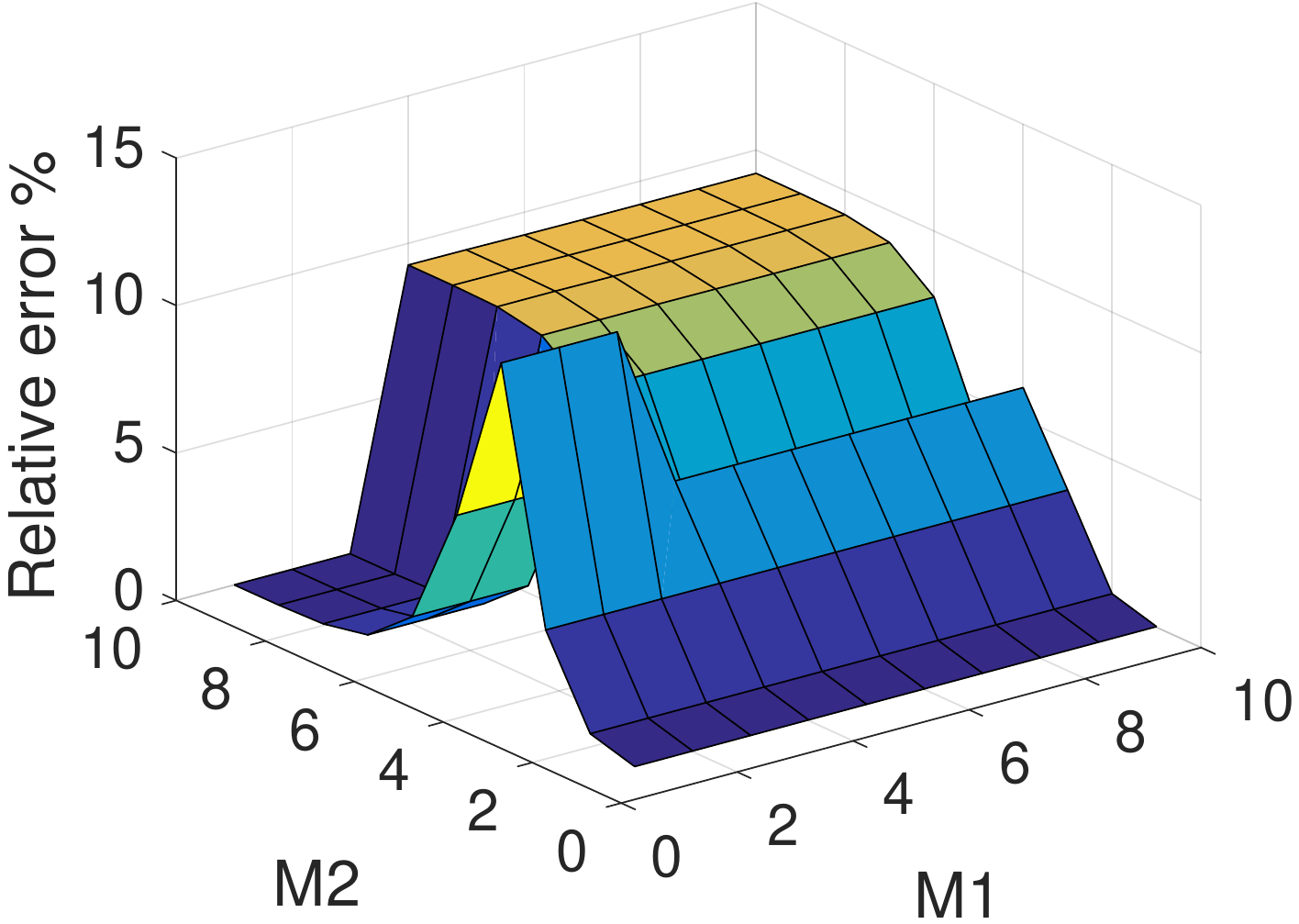}
	\caption{Distflow in case $\blt{K}=(10,10)$ and $\blt{\lambda}=(12,12)$}
\label{fig1}
\end{figure}

\begin{figure}[!h]
	\centering
		\includegraphics[width=.45\linewidth]
		{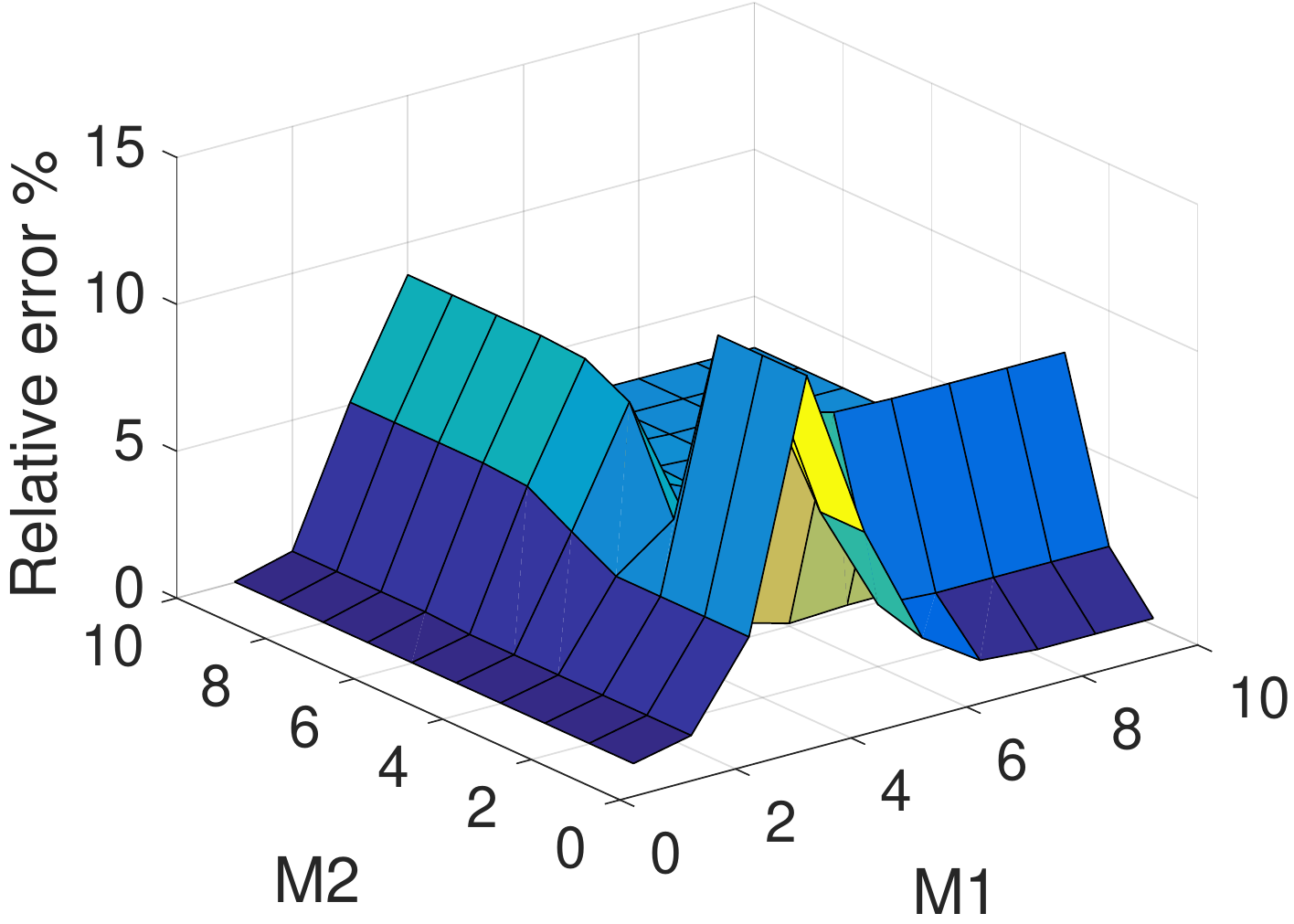}
		\includegraphics[width=.45\linewidth]
		{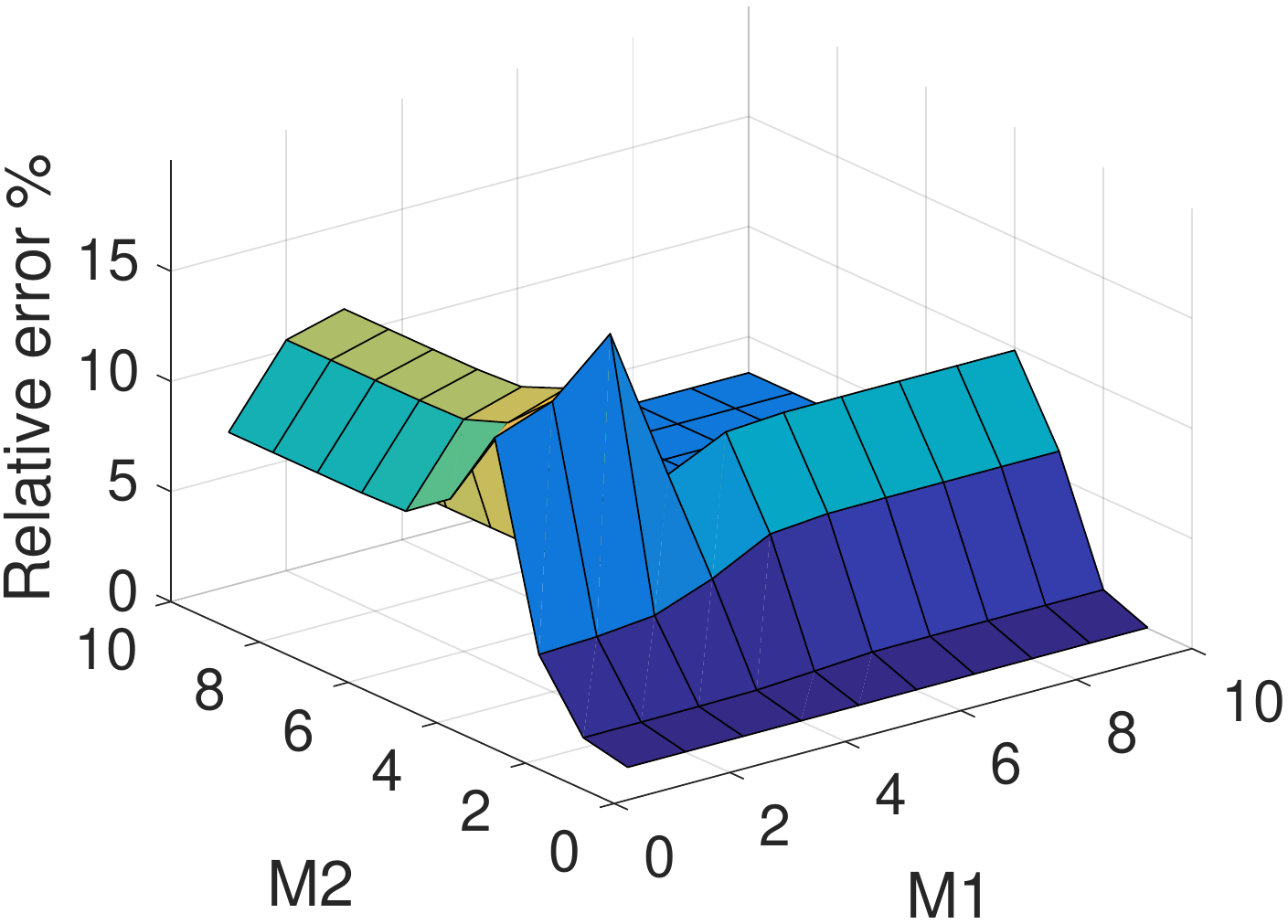}
	\caption{AC in case $\blt{K}=(10,10)$ and $\blt{\lambda}=(12,12)$}
\label{fig2}
\end{figure}

Now, we move to a model with all the additional constraints. In this case, we fix $c^{\text{max}}=1$ and we plot the relative error between the fluid approximation and the stochastic model for the fraction of EVs that get successfully charged. Using arguments from queueing theory, it can be shown that this fraction equals $1-\E{Z_i}/\E{Q_i}$. For the fluid approximation, we replace the numerator with $z_i^*$, while $\E{Q_i}$ can be computed
explicitly, using the Erlang Loss formula \cite{kelly2014stochastic}. Figures~\ref{fig1} and \ref{fig2} show the results for all possible values of $M_i$ and for both load models. Though the quality of the fluid approximation deteriorates, the relative error is generally below  $10\%$ and for reasonably high values of $M_i$, it is even smaller. We also expect these results to improve for bigger $K$.
For  higher values of $ K$ one needs to solve millions of OPF problems of type \eqref{OP} to obtain the steady-state behavior of the stochastic model, while the fluid model only requires the solution of a single OPF problem.


\section{Proofs}\label{proofs}
\section*{Proof of Theorem~\ref{Th:UnFP}}\label{sec:PrTheorem1}
Define $a_{ij}=\inf\{x\geq 0: g_{ij}(x) =\gamma_{ij}/\E{D_j}\}$. By \cite[Lemma 6]{remerova14}, we have that $g_{ij}(\cdot)$ is continuous, strictly increasing in $[0,a_{ij}]$, and constant in $(a_{ij}, \infty]$. Further, by \cite[Lemma 7]{remerova14}, we have that if $a_{ij}<\infty$, then $\inf{(D_j/B_{j})}= 1/a_{ij}$ and if $a_{ij}=0$, then $\inf{(D_j/B_{j})}= 0$.

 For $i\in \Node \setminus \{0\}$, we have that
\begin{equation}\label{pr:help}
 \frac{\partial }{\partial p_{ij}}W_{kk}(\blt {z}, \blt p)=z_{ij}\frac{\partial }{\partial \Lambda_{ij}}W_{kk}(\blt {\Lambda}).
\end{equation}
 To see this, observe that by \eqref{eq:AC} and the definition of the bandwidth allocation $\blt \Lambda$, as $V_k(\blt {z}, \blt {p})$ depends on the vectors $\blt{z}$ and $\blt{p}$ only through the product $z_{ij}p_{ij}$, we get
$$V_k(\blt z, \blt p)=V_k(\blt \Lambda).$$
Applying the transformation
$\Lambda_{ij}(\blt z,\blt p)=z_{ij} p_{ij}$, we have
\begin{equation*}
\frac{\partial }{\partial p_{ij}}
V_k(\blt {z},\blt p)
=\frac{\partial }{\partial \Lambda_{ij}} V_k(\blt {\Lambda})
\frac{\partial }{\partial p_{ij}} \Lambda_{ij}
=z_{ij} \frac{\partial }{\partial \Lambda_{ij}}V_k(\blt {\Lambda}).
\end{equation*}
Recalling that $W_{kk}=V_kV_k$, we obtain that
\begin{equation*}
\begin{split}
  \frac{\partial }{\partial p_{ij}}W_{kk}(\blt{z},\blt{p}) =&
  2 V_k(\blt {z},\blt p)
  \frac{\partial }{\partial p_{ij}}V_k(\blt {z},\blt p)\\
=&2z_{ij}V_k(\blt {\Lambda})
  \frac{\partial }{\partial \Lambda_{ij}}V_k(\blt {\Lambda})
  =z_{ij}\frac{\partial }{\partial \Lambda_{ij}}W_{kk}(\blt{\Lambda}).
\end{split}
\end{equation*}

By the KKT conditions for \eqref{OP} there exist Lagrangian multipliers  $(h_{i}^1, h_i^2, h_{ij}^3, \bar h_{ij}^4)
\in \R^{2I+ 2(I\times J)}_{+}$
such that for $i\in\Node \setminus \{0\}$, $j\in J$,
\begin{equation*}
\begin{split}
z_{ij}^*  \frac{\partial }{\partial p_{ij}}u_{ij}(p_{ij}(\blt z^*))=
 \sum_{k=1}^{I}
  \frac{\partial }{\partial p_{ij}}W_{kk}(\blt p(\blt z^*))
 (h_{k}^1-
 h_{k}^2)\\
  +h_{ij}^3 z_{ij}^*+ \bar h_{ij}^4,
 \end{split}
 \end{equation*}
and
\begin{align*}
  h_{i}^1 (W_{ii}(\blt p(\blt z^*))-\overline{\upsilon}_i) &=0, \quad
  h_{i}^2 (W_{ii}(\blt p(\blt z^*))-\underline{\upsilon}_i)=0,  \\
  h_{ij}^3 (z_{ij}^*p_{ij}(\blt z^*)-M_i)&=0, \quad
  \bar h_{ij}^4 (p_{ij}(\blt z^*)-c_j^{max}) =0.
\end{align*}
Setting $h_{ij}^4= \bar h_{ij}^4/z_{ij}^*$, the previous equations  take the following (equivalent) form
\begin{equation*}
\begin{split}
\frac{\partial }{\partial p_{ij}}u_{ij}(p_{ij}(\blt z^*))  =
 \frac{1 }{z_{ij}^*}\sum_{k=1}^{I} \frac{\partial }{\partial p_{ij}}W_{kk}(\blt p(\blt z^*))
 (h_{k}^1-
 h_{k}^2)\\
 +h_{ij}^3 + h_{ij}^4,
 \end{split}
\end{equation*}
and
\begin{align*}
  h_{i}^1 (W_{ii}(\blt p(\blt z^*))-\overline{\upsilon}_i) &=0, \quad
  h_{i}^2 (W_{ii}(\blt p(\blt z^*))-\underline{\upsilon}_i))=0,  \\
  h_{ij}^3 (z_{ij}^*p_{ij}(\blt z^*)-M_i)&=0, \quad
  h_{ij}^4 (p_{ij}(\blt z^*)-c_j^{max}) =0.
\end{align*}
By the definition of $g_{ij}(\cdot)$ and \eqref{eq:fluid proxy}, we have that $\Lambda_{ij}(z^*)=g_{ij}(p_{ij}(\blt z^*))$. Moreover, by the assumption in Theorem~\ref{Th:UnFP}, we have that $1/a_{ij}=\inf{(D_j/B_{j})}\leq 1/c_j^{\text{max}}$. That is, $c_j^{\text{max}}\leq a_{ij}$.
Thus, $g_{ij}(\cdot)$ is strictly increasing in $[0,c_j^{\text{max}}]$. This implies
$p_{ij}(\blt z^*)=g_{ij}^{-1}(\Lambda_{ij}(z^*))$. We note that
$(p_{ij}(\blt z^*)-c_j^{max}) =0$ if and only if
$(\Lambda_{ij}(\blt z^*)-g_{ij}(c_i^{max}))=0$.
Using the last observations and \eqref{pr:help}, the above equations can be rewritten as follows
\begin{equation*}
\begin{split}
 \frac{\partial }{\partial \Lambda_{ij}}u_{ij}(g_{ij}^{-1}(\Lambda_{ij}(\blt z^*)))  =
\sum_{k=1}^{I} \frac{\partial }{\partial \Lambda_{ij}}W_{kk}(\blt \Lambda)
 (h_{k}^1-
 h_{k}^2)\\
 +h_{ij}^3 + h_{ij}^4,
 \end{split}
\end{equation*}
and
\begin{align*}
  h_i^1 (W_{ii}(\blt \Lambda(\blt z^*))-\overline{\upsilon}_i) &=0, \quad
  h_i^2 (W_{ii}(\blt \Lambda(\blt z^*))-\underline{\upsilon}_i)=0, \\
  h_{ij}^3 (\Lambda_{ij}(\blt z^*)-M_i)&=0,  \quad
   h_{ij}^4 (\Lambda_{ij}(\blt z^*)-g_{ij}(c_i^{max}))=0.
\end{align*}
Now, we observe that the last equations are the KKT conditions for~\eqref{GDOP}.
To complete the proof, it remains to be shown that the function $G_{ij}(\cdot)$ is strictly concave. To this end, observe that $g_{ij}^{-1}(\cdot)$ is strictly increasing and $ u_{ij}'(\cdot)$ is strictly decreasing since $ u_{ij}(\cdot)$ is a strictly concave function. It follows that $G_{ij}'(\cdot)$ is strictly decreasing and hence $G_{ij}(\cdot)$ is  a strictly concave function. The latter implies that the optimization problem \eqref{GDOP} has a unique solution $\blt{\Lambda}^*$, which is independent of the invariant point $\blt{z}^*$.
Further, the unique invariant point is given by
\begin{equation*}\label{eq: InP}
  z_{ij}^*=\frac{\Lambda_{ij}^*}{p_{ij}(\blt z^*)}
  =\frac{\Lambda_{ij}^*}{g_{ij}^{-1}(\Lambda_{ij}^*)}.
\end{equation*}

\section*{Proof of Proposition~\ref{eq:SSL}}\label{sec:eq:SSL}
Given the similarities with \cite[Theorem~3, Corollary~1]{remerova14}, we just sketch the proof. Under the assumption that
 $K_i=\infty$, we have that
$\ind {q_i(t)<K_i}=1$.
 By following the same argument as in \cite[Theorem~3]{remerova14}, there exist $\blt l, \blt h \in (0,\infty)^{I\times J}$, such that
\begin{equation}\label{in:san}
  l_{ij}\leq \liminf_{t\rightarrow \infty} z_{ij}(t)
  \leq \limsup_{t\rightarrow \infty} z_{ij}(t) \leq h_{ij}.
\end{equation}
By monotonicity, $\blt l$ and $\blt h$ satisfy
\begin{align*}
l_{ij}=\gamma_{ij}\E{ \min \{D_j , \frac{ B_{j}} {p_{ij(\blt l)}}\}},\
h_{ij}=\gamma_{ij}\E{ \min \{D_j , \frac{ B_{j}} {p_{ij(\blt h)}}\}}.
\end{align*}
By the uniqueness of the invariant point, we have that $l_{ij}=h_{ij}$. Now, by \eqref{in:san}, the result follows.

\section*{Proof of Proposition~\ref{prop:product form}}
Using the KKT conditions,
it can be shown that
$p_i(\blt z)= \frac{\delta}{\SR_i \sum_{l=1}^{I}z_l}.$
Next, observe that the so-called balance property \cite{bonald2002Insensitivity} holds: for
$i,k\in \Node \setminus\{0\}$,
\begin{equation*}
p_i(\blt z+\blt{e}_k)p_k(\blt z)=p_i(\blt z) p_k(\blt z+\blt{e}_i),
\end{equation*}
where $\blt e_i$, $i=1,\ldots,I$, denote unit vectors.

If $B$ is an exponential random variable, the process $\blt Z(t),\ t\geq 0$ is Markov.
Let $\blt Z(t)=\blt z$  denote the state of the process $\blt Z(t)$ at time $t\geq0$. The process can move to state $\blt z+\blt e_i$ with rate $\lambda_i$ and to state $\blt z-\blt e_i$ with rate $z_i p_i(\blt z)/\E{B}$, if $z_i>0$. Setting $\mu_i=\frac{\E{B}\SR_i}{\delta}$, the transition rates become the same as those of a multiclass Markovian processor-sharing queue with $I$ classes of customers and mean service times $\mu_i$. By \cite{cohen79}, the stationary distribution of a multiclass processor-sharing queue is given by \eqref{eq:PS}.

Further, it is shown in \cite{bonald2002Insensitivity}, that the stationary distribution of a multiclass processor-sharing queue is insensitive to the distribution of service times, if the balance property is satisfied. That is, \eqref{eq:PS} holds for general charging requirements $B$. The construction in \cite{bonald2002Insensitivity} can also be carried out in our setting.

\section*{Proof of Proposition~\ref{prop:explicitZ}}
For simplicity take $\E{B}=1/\mu$ and $\E{D}=1/\nu$, and $z_i(0)=0$. By the assumption that $K_i=\infty$, it follows that
 $\gamma_i(t)\equiv \gamma_i\equiv \lambda_i$ (for the definition of $\gamma_i$, see \eqref{eq:defGamma}). Before we move to the main part of the proof, we show some helpful relations.

Under the assumptions in Section~\ref{Additional examples under a line topology}, we can solve \eqref{GDOP} explicitly. To see this, note that
$G_i'(\Lambda_i)=w_i\frac{\gamma_i-\mu \Lambda_i}{\nu \Lambda_i}$.
Writing the KKT conditions for~\eqref{GDOP}, we have that there exists an $h\in R_+$, such that for any $i \in \Node\setminus \{0\}$, the optimal solution $\blt \Lambda^*$ satisfies
$w_i\frac{\gamma_i-\mu \Lambda_i^*}{\nu \Lambda_i^*}=h\SR_i $,
and
\begin{equation}\label{eq:FC}
\sum_{k=1}^{I}\SR_k \Lambda_k^*=\delta.
\end{equation}
The first equation (taking into account $w_i=\SR_i$) implies that
$\frac{\gamma_i-\mu \Lambda_i^*}{\nu \Lambda_i^*}=
\frac{\gamma_1-\mu \Lambda_1^*}{\nu \Lambda_1^*}$.
Finally, we get the following expressions for $\Lambda_i^*$,
\begin{equation}\label{eq:ExpL}
\Lambda_i^*=\frac{\gamma_i}{\gamma_1}\Lambda_1^*,
\end{equation}
and by \eqref{eq:FC},
\begin{equation}\label{eq:FA}
\Lambda_i^*= \frac{\gamma_i\delta}{\sum_{k=1}^{I}\SR_k\gamma_k}.
\end{equation}

Next, by \eqref{eq:fluid proxy}, we compute the invariant point,
which is given by
\begin{equation}\label{eq:FPExplicit}
  z_i^*=\frac{\gamma_i-\mu \Lambda_i^*}{\nu}.
\end{equation}
It is helpful to note the following relation for the invariant point. Combining the last equation and \eqref{eq:ExpL}, we get
\begin{equation}\label{eq:FP}
  z_i^*=\frac{\gamma_i-\mu (\gamma_i/\gamma_1)\Lambda_1^*}{\nu}
  =\frac{\gamma_i}{\gamma_1}\frac{\gamma_1-\mu \Lambda_1^*}{\nu}
  =\frac{\gamma_i}{\gamma_1}z_1^*.
\end{equation}

Now, we move to the main part of the proof.
Under the Markovian assumptions, $z_i(\cdot)$ is given (alternatively) by the following ODE:
\begin{equation*}
z_i'(t)=\gamma_i-\nu z_i(t)-\mu z_i(t)\frac{\delta}{\sum_{k=1}^{I}\SR_k z_k(t)}.
\end{equation*}
The last ODE has a unique solution $\blt {z}(\cdot)$ for given initial point. So, it is enough to show that the function $z_i(t)=z_i^*(1-e^{-\nu t})$ satisfies the previous ODE. Plugging in it into the ODE, we have that
\begin{equation*}
\nu z_i^*e^{-\nu t}=\gamma_i-\nu z_i^*(1-e^{-\nu t})
-\frac{\mu z_i^*(1-e^{-\nu t})\delta}{\sum_{k=1}^{I}\SR_k z_k^*(1-e^{-\nu t})},
\end{equation*}
which can be simplified to
$\gamma_i=\nu z_i^*
+\frac{\mu z_i^*\delta}
{\sum_{k=1}^{I}\SR_k z_k^*}$.
By \eqref{eq:FPExplicit}, we derive
$\Lambda_i^*=\frac{ z_i^*\delta}
{\sum_{k=1}^{I}\SR_k z_k^*}$.
Now, we apply \eqref{eq:FP} to get
$\Lambda_i^*=\frac{ \gamma_i \delta}
{\sum_{k=1}^{I}\SR_k \gamma_k}$,
which holds by \eqref{eq:FA}. 

\section*{Proof of Proposition~\ref{prop:load balancing}}
Under the assumptions of Proposition~\ref{prop:load balancing}, and observing that $W^{lin}_{ii}$ is decreasing in $i$, \eqref{OP} takes the following form:
\begin{equation*}
\begin{aligned}
& \underset{}{\text{max}_{\blt p}}
& & \sum_{i=1}^{I} \sum_{j=1}^{J} z_{ij} w_{ij} u(p_{ij}) \\
& \text{subject to}
&& \sum_{k=1}^{I}R_{k-1k}\sum_{m=k}^{I}
  \sum_{j=1}^{J} z_{mj}p_{mj} \leq \delta .\\
\end{aligned}
\end{equation*}
By the Karush-Kuhn-Tucker (KKT) conditions, there exists an
$h\in \R_{+}$ such that for any $i\in \Node \setminus \{0\}$ and $j\in \Edge$,
\begin{equation}\label{eq:LP1}
z_{ij}w_{ij}u'(p_{ij})= h z_{ij} \SR_i,
\end{equation}
and
\begin{equation}\label{eq:ConKKT}
h\Big(\sum_{k=1}^{I}R_{k-1k}\sum_{m=k}^{I}
  \sum_{j=1}^{J} z_{mj}p_{mj}-\delta\Big)=0.
\end{equation}
Note that $h$ cannot be zero due to \eqref{eq:LP1}. Again, by \eqref{eq:LP1}, we have that
$\frac{w_{ij}u'(p_{ij})}{ \SR_i}$ should be constant for any $i,j$. In particular,
\begin{equation}\label{eq:EqL}
\frac{w_{ij}u'(p_{ij})}{ \SR_i}=\frac{w_{11}u'(p_{11})}{ \SR_1}.
\end{equation}
Choosing $w_{ij}=\SR_i$, and noting that $u'(\cdot)$ is a strictly decreasing function, we have that $p_{ij}=p_{11}$, for all $i,j \geq 1$.
Combining \eqref{eq:ConKKT} and \eqref{eq:EqL}, we derive the expression for $p(\blt z)$.

\ifCLASSOPTIONcaptionsoff
  \newpage
\fi

\begin{IEEEbiography}
[{\includegraphics[width=1in,height=1.25in,clip,keepaspectratio]{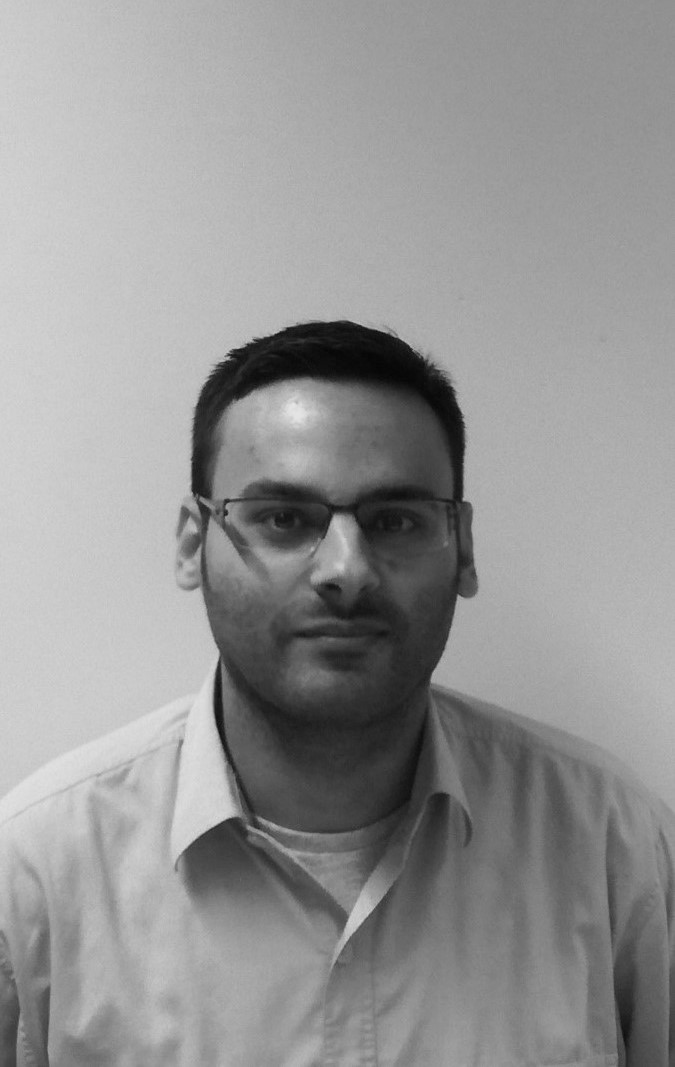}}]
{Angelos Aveklouris}
studied Mathematics at the National and Kapodistrian University of Athens,
where he graduated in 2012. In 2015, he earned a masters degree in Mathematical Modeling in Modern Technologies and Economics
from the department of Mathematics at the National Technical University of Athens.
Later that year, he began his PhD research project in the department of Mathematics and Computer Science at the Eindhoven University of Technology.
His research interests are mainly in stochastic networks and probability theory.
\end{IEEEbiography}

\begin{IEEEbiography}
[{\includegraphics[width=1in,height=1.25in,clip,keepaspectratio]{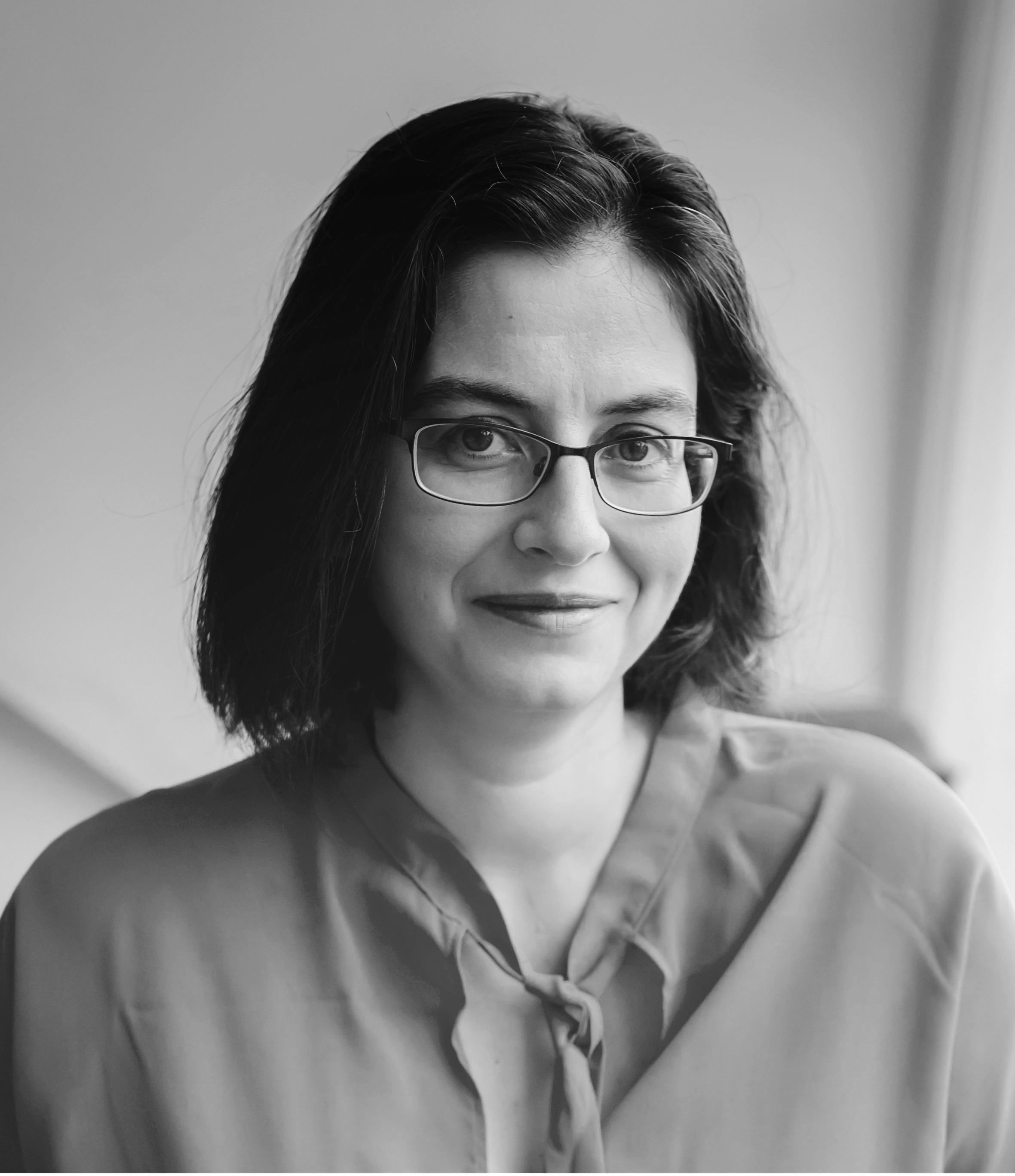}}]
{Maria Vlasiou}
 (BSc Mathematics 2002, Hons.; PhD Mathematics 2006) is an Associate Professor at the Eindhoven University of Technology, a scientific staff member at CWI Amsterdam, and a research fellow of EURANDOM. Her research focuses on the performance of stochastic processing networks with layered architectures and on perturbation analysis for heavy-tailed risk models. Dr. Vlasiou serves as the President of the Dutch association of female professional mathematicians (EWM-NL), as a member of the scientific advisory board in mathematics of the Lorenz Center, and as associate editor for IISE Transactions and for Mathematical Methods of Operations Research.

\end{IEEEbiography}
\begin{IEEEbiography}
[{\includegraphics[width=1in,height=1.25in,clip,keepaspectratio]{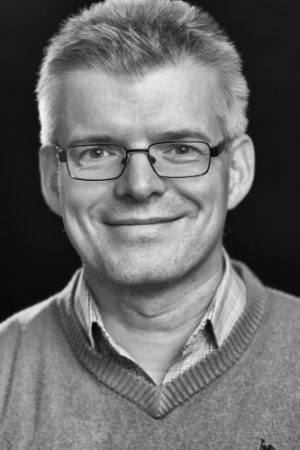}}]{Bert Zwart}
(MA Econometrics 1997, PhD Mathematics 2001, Member
of IEEE since 2017) is leader of the CWI Stochastics group (Amsterdam) and Professor at the Eindhoven University of Technology. He is Stochastic Models area editor for Operations Research, the flagship journal of his profession, since 2009. His research expertise is in applied probability and stochastic networks. His work on power systems is focusing on the applications of probabilistic methods to reliability issues and rare events such as cascading failures and blackouts. He is co-organizer of a special semester on the Mathematics of energy systems taking place in Cambridge UK, spring 2019.
\end{IEEEbiography}
\end{document}